\documentclass{amsart}
\usepackage{indentfirst}
\usepackage{bm}
\usepackage{amssymb}
\usepackage{stmaryrd}
\usepackage{enumerate}
\usepackage{amsmath}
\usepackage{commath}
\usepackage{hyperref}
\usepackage{cleveref}
\usepackage{graphicx}
\usepackage{amsaddr}
\usepackage[caption=false]{subfig}
\usepackage{xspace,color}

\newcommand{\jump}[1]{\llbracket #1 \rrbracket}

\makeatletter
\newcommand{\tnorm}{\@ifstar\@tnorms\@tnorm}
\newcommand{\@tnorm}[2][]{%
	\mathopen{#1|\mkern-1.5mu#1|\mkern-1.5mu#1|}
	#2
	\mathclose{#1|\mkern-1.5mu#1|\mkern-1.5mu#1|}
}
\newtheorem{remark}{Remark}
\newtheorem{theorem}{Theorem}
\newtheorem{lemma}{Lemma}
\begin{document}
\title[Analysis of EDG-HDG for Biot]{Analysis of an embedded-hybridizable discontinuous Galerkin method for Biot's consolidation model}

\author{Aycil Cesmelioglu$^1$, Jeonghun J. Lee$^2$, Sander Rhebergen$^3$ }
\address{$^1$Department of Mathematics and Statistics, Oakland
  University, Michigan, USA \\ 
  $^2$ Department of Mathematics, Baylor University, Waco, Texas, USA \\
  $^3$ Department of Applied Mathematics, University of
  Waterloo, Ontario, Canada }

\email{$^1$cesmelio@oakland.edu, $^2$jeonghun\_lee@baylor.edu, $^3$srheberg@uwaterloo.ca}

\begin{abstract}
  We present an embedded-hybridizable discontinuous Galerkin finite
  element method for the total pressure formulation of the
  quasi-static poroelasticity model. Although the displacement and the
  Darcy velocity are approximated by discontinuous piece-wise
  polynomials, $H(\text{div})$-conformity of these unknowns is
  enforced by Lagrange multipliers. The semi-discrete problem is shown
  to be stable and the fully discrete problem is shown to be
  well-posed. Additionally, space-time a priori error estimates are
  derived, and confirmed by numerical examples, that show that the
  proposed discretization is free of volumetric locking.
\end{abstract}

\date{July, 2022}
\maketitle

\maketitle

\section{Introduction}
\label{sec:introduction}

Poroelasticity models are systems of partial differential equation that describe the physics of deformable porous media saturated by fluids. They were originally developed for geophysics applications in petroleum engineering but nowadays they are also widely used for biomechanical modeling. The first poroelasticity models were derived by Biot \cite{Biot-low-frequency,Biot-high-frequency}. Since then, mathematical properties and numerical methods for these models have been widely studied. Here we give a brief literature review.

Early studies on linear poroelasticity models include well-posedness analysis and finite element discretizations for quasi-static \cite{Zenisek:1984,Showalter:2000} and dynamic \cite{Santos:1986a,Santos:1986b,Zienkiewicz-Shiomi:1984} models. For quasi-static models with incompressible elastic grains, Murad et al. \cite{Murad:1996} observed spurious pressure oscillations of certain finite element discretizations for small time and studied their asymptotic behavior. Phillips and Wheeler \cite{Phillips:2009} connected these pressure oscillations to volumetric locking due to incompressibility of the displacement. They further developed numerical methods in \cite{Phillips:2008,Phillips:2009} coupling mixed methods and discontinuous Galerkin methods that do not show pressure oscillations. Yi \cite{Yi:2013,Yi:2014,Yi:2017} proposed numerical methods coupling mixed and nonconforming finite elements that are also free of pressure oscillations. An analysis to address the volumetric locking problem for poroelasticity was first presented in \cite{JJLee:2016} adopting mixed methods for linear elasticity. Various numerical methods avoiding this locking problem have since been studied using nonconforming or stabilized finite elements \cite{Lee:2017,Rodrigo:2016,Hu:2017,Boffi:2016}, the total pressure formulation \cite{JJLee:2017,Oyarzua:2016,Feng:2018}, and exactly divergence-free finite element spaces \cite{Kanschat:2018,Hong:2018}. A non-symmetric interior penalty discontinuous Galerkin method was numerically shown to be locking free for small enough penalty parameter in \cite{Riviere:2017}.

Discontinuous Galerkin methods are known to be computationally expensive. A remedy for this was provided by Cockburn et el. \cite{Cockburn:2009a} by introducing the hybridizable discontinuous Galerkin (HDG) framework for elliptic problems. Indeed, element unknowns can be eliminated from the problem resulting in a global problem for facet unknowns only. The number of globally coupled degrees-of-freedom can be reduced even further using the embedded discontinuous Galerkin (EDG) framework \cite{Cockburn:2009,Guzey:2007}; where the HDG method uses a discontinuous trace approximation, the EDG method uses a continuous trace approximation. HDG, and related hybrid high-order (HHO), methods have recently been introduced for the poroelasticity problem \cite{Fu:2018-biot,Kraus:2021,Botti:2021}. These discretizations consider the primal bilinear form for linear elasticity. In contrast, in this paper we adopt the total pressure formulation  \cite{JJLee:2017,Oyarzua:2016} and present novel HDG and EDG-HDG methods for the quasi-static poroelasticity models. (It is possible to also consider an EDG method for the poroelasticity model, however, such a discretization is sub-optimal.) The total pressure formulation provides a natural decoupling of the linear elasticity and Darcy equations in the incompressible limit. Indeed, in this limit our discretizations reduce to the exactly divergence-free HDG and EDG-HDG discretizations of \cite{Rhebergen:2017,Rhebergen:2020} for the Stokes problem and the hybridized formulation of \cite{Arnold:1985} for the Darcy problem. We further remark that the total pressure formulation has been applied also in the context of magma/mantle dynamics problems \cite{Katz:2007,Keller:2013} where it was shown to be advantageous in the context of coupled physics problems beyond quasi-static poroelasticity problems. 

We present an analysis of the proposed HDG and EDG-HDG methods in which we show that the space-time discretizations are well-posed. We further determine an a priori error estimate for all unknowns that is robust in the incompressible limit and for arbitrarily small specific storage coefficient. We remark that the standard approach of analyzing time-dependent problems is to use discrete Gr\"onwall inequalities. However, this results in error bounds with a coefficient that grows exponentially in time. We present an alternative approach that avoids this exponential term. 

The remainder of this paper is organized as follows. We present Biot's consolidation model in \cref{sec:biot}. The HDG and EDG-HDG methods for Biot's model is presented in \cref{sec:hdg} together with a stability proof for the semi-discrete problem. Well-posedness and a priori error estimates for the fully discrete problem are shown in \cref{sec:apriori}. The analysis is verified by numerical examples in \cref{sec:numerical_examples} and conclusions are drawn in \cref{sec:conclusions}.

\section{Biot's consolidation model}
\label{sec:biot}
To introduce Biot's consolidation model, let us introduce the
following notation. Let $\Omega \subset \mathbb{R}^d$, $d=2,3$ be a
bounded polygonal domain with a boundary partitioned 
as $\partial \Omega = \overline{\Gamma}_P \cup \overline{\Gamma}_F$ and
$\partial \Omega = \overline{\Gamma}_D \cup \overline{\Gamma}_T$,
where $\Gamma_P \cap \Gamma_F = \emptyset$, $|\Gamma_P| > 0$,
$\Gamma_D \cap \Gamma_T = \emptyset$, and $|\Gamma_D| > 0$. We denote
the unit outward normal to $\partial \Omega$ by $n$ and we denote by
$I = (0, T]$ the time interval of interest.

Let $f:\Omega\times I \to \mathbb{R}^d$ be a given body force and let
$g:\Omega\times I \to \mathbb{R}$ be a given source/sink
term. Furthermore, let $\kappa>0$ be a scalar constant that represents
the permeability of the porous media, $c_0 \ge 0$ the specific storage
coefficient, and $0 < \alpha < 1$ the Biot--Willis constant. Denoting
Young's modulus of elasticity by $E$ and Poisson's ratio by $\nu$, in
the case of plane strain, the Lam\'e constants are given by
$\lambda = E\nu/((1+\nu)(1-2\nu))$ and $\mu = E/(2(1+\nu))$.

Biot's consolidation model describes a system of equations for the
displacement of the porous media, $u:\Omega\times I\to\mathbb{R}^d$,
and the pore pressure of the fluid $p:\Omega\times I\to
\mathbb{R}$. Denoting by
$\sigma = 2\mu\varepsilon(u) + \lambda\nabla\cdot u\mathbb{I} - \alpha
p \mathbb{I}$ the total Cauchy stress, where $\mathbb{I}$ is the $d\times d$-dimensional identity matrix, this model is given by
\begin{equation}
  \label{eq:originalbiotmodel}
  -\nabla\cdot \sigma = f,
  \qquad
  \partial_t(c_0 p + \alpha\nabla\cdot u) - \nabla\cdot(\kappa\nabla p) = g,
  \qquad
  \text{ in } \Omega \times I.
\end{equation}
Following \cite{JJLee:2017}, by introducing the total pressure
$p_T:= -\lambda \nabla \cdot u+\alpha p$ and the Darcy velocity
$z := -\kappa\nabla p$, we may write Biot's consolidation model also
as:
\begin{subequations}
  \label{eq:biot}
  \begin{align}
    \label{eq:biot_a}
    -\nabla\cdot 2\mu\varepsilon(u) + \nabla p_T
    &= f & & \text{in}\ \Omega \times I,
    \\
    \label{eq:biot_b}
    - \nabla\cdot u - \lambda^{-1} (p_T-\alpha p) 	&= 0
                                             & & \text{in}\ \Omega \times I,
    \\
    \label{eq:biot_c}
    \partial_t\del[1]{ c_0p + \lambda^{-1}\alpha (\alpha p - p_T)} + \nabla\cdot z &= g
                                & & \text{in}\ \Omega \times I,
    \\
    \label{eq:biot_d}
    \kappa^{-1} z + \nabla p &= 0
                                             & & \text{in}\ \Omega \times I,
  \end{align}
\end{subequations}
which will be the formulation studied in this article. Noting that
$\sigma = 2\mu\varepsilon(u) - p_T\mathbb{I}$, we close the model by
imposing the following boundary and initial conditions:
\begin{subequations}
  \begin{align}
    \label{eq:bics-u}
    u &= 0 && \text{on}\  \Gamma_D\times I,
    \\
    \label{eq:bics-p}
    p &= 0 && \text{on}\ \Gamma_P\times I,
    \\
    \label{eq:bics-z}
    z\cdot n &= 0 && \text{on}\ \Gamma_F\times I,
    \\
    \label{eq:bics-sigma}
    \sigma n &= 0 &&  \text{on}\ \Gamma_T\times I,
    \\
    \label{eq:ic_p}
    p(x,0)& =p_0(x) && \text{in}\ \Omega,
    \\
    \label{eq:ic_u}
    u(x,0) &= u_0(x) && \text{in}\ \Omega.
  \end{align}
  \label{eq:bics}
\end{subequations}

In the remainder of this article we assume that $c_0$, $\mu^{-1}$,
$\kappa$, and $\mu$ are bounded above by a constant $C$. We
furthermore assume that there exists a $\nu_*$ such that
$0 < \nu_* \le \nu < 0.5$ on $\Omega$. As a consequence,
$C_* \mu \le \lambda$ with $C_* = 2\nu^*/(1-2\nu^*)$.

\section{The embedded-hybridizable discontinuous Galerkin method}
\label{sec:hdg}

\subsection{Notation}
On a Lipschitz domain $D$ in $\mathbb{R}^d$, we denote by $W^{l,p}(D)$
the usual Sobolev spaces for $l \ge 0$ and $1 \le p \le \infty$ (see,
for example, \cite{Adams:book}). When $p=2$, we define on
$H^l(D) = W^{l,2}(D)$ the norm $\norm[0]{\cdot}_{l, D}$ and semi-norm
$\envert[0]{\cdot}_{l, D}$. We note that $L^2(D) = H^0(D)$ is the
Lebesque space of square integrable functions with norm
$\norm[0]{\cdot}_D = \norm[0]{\cdot}_{0,D}$ and inner product
$(\cdot, \cdot)_D$. Vector-valued function spaces will be denoted by
$\sbr[0]{L^2(D)}^d$ and $\sbr[0]{H^l(D)}^d$. The $L^2$-inner product
over a surface $S \subset \mathbb{R}^{d-1}$ will be denoted by
$\langle \cdot, \cdot \rangle_S$.

Let $X$ be a Banach space and $J=(0,T]$, $T>0$ a time interval. We
denote by $C^0(J; X)$ the space of continuous functions $f:J \to X$,
which is equipped with the norm
$\norm{f}_{C^0(\overline{J};X)} := \sup_{t\in \overline{J}}
\norm[0]{f(t)}_X$. By $C^k(J; X)$, $k \ge 0$, we denote the space of
continuous functions $f:J \to X$ such that
$\partial_t^i f \in C^0(J, X)$ for $1 \le i \le k$. For
$1 \le p < \infty$, $W^{k,p}(J; X)$ is defined to be the closure of
$C^k(J; X)$ with respect to the norm
\begin{equation*}
  \norm[0]{ f }_{W^{k,p}(J; X)}^p := \int_0^T \sum_{i=0}^k \norm[0]{ \partial_t^i f(t) }_X^p \dif t.
\end{equation*}
We note that for $k=0$, $W^{k,p}(J; X) = L^p(J; X)$.

Let $\mathcal{T}_h$ be a family of shape-regular simplicial
triangulations of the domain $\Omega$. We will denote the diameter of
an element $K \in \mathcal{T}_h$ by $h_K$, the meshsize by
$h := \max_{K \in \mathcal{T}_h} h_K$, and the sets of interior facets
and facets that lie on $\Gamma_D$, $\Gamma_P$, $\Gamma_F$, and
$\Gamma_T$ by, respectively, $\mathcal{F}_h^i$, $\mathcal{F}_h^D$,
$\mathcal{F}_h^P$, $\mathcal{F}_h^F$, and $\mathcal{F}_h^T$. The set
of all facets is denoted by $\mathcal{F}_h$ and their union is denoted by $\Gamma_0$. On the boundary of an element $K$, we denote by $n_K$ the outward unit normal vector, although, where no confusion will occur we drop the subscript $K$. On the mesh and skeleton we define the inner products
\begin{align*}
  (\phi, \psi)_{\Omega}
  &:=\sum_{K\in \mathcal{T}_h}(\phi, \psi)_K,
  & \langle\phi, \psi\rangle_{\partial \mathcal{T}_h}
  &:=\sum_{K\in \mathcal{T}_h}\langle \phi, \psi \rangle_{\partial K},
  && \text{ if } \phi,\psi \text{ are scalar},
  \\
  (\phi, \psi)_{\Omega}
  &:=\sum_{i=1}^d(\phi_i, \psi_i)_{\Omega},
  & \langle\phi, \psi\rangle_{\partial \mathcal{T}_h}
  &:=\sum_{i=1}^d\langle\phi_i, \psi_i\rangle_{\partial \mathcal{T}_h},
  && \text{ if } \phi,\psi \text{ are vector-valued}.
\end{align*}
The norms induced by these inner products are denoted by
$\norm[0]{\cdot}_{\Omega}$ and $\norm[0]{\cdot}_{\partial \mathcal{T}_h}$,
respectively.

Sets of polynomials of degree not larger than $l \ge 0$ defined on,
respectively, an element $K \in \mathcal{T}_h$ and a facet
$F \in \mathcal{F}_h$ will be denoted by $\mathbb{P}_l(K)$ and
$\mathbb{P}_l(F)$. As approximation spaces we then use:
\begin{equation}
  \label{eq:approximationspaces}
  \begin{split}
    V_h
    &:= \cbr[1]{v_h\in \sbr[0]{L^2(\Omega)}^d
      : \ v_h \in \sbr[0]{\mathbb{P}_k(K)}^d, \ \forall\ K\in\mathcal{T}_h},
    \\
    \bar{V}_h
    &:= \cbr[1]{\bar{v}_h \in \sbr[0]{L^2(\Gamma_0}^d:\ \bar{v}_h \in
      \sbr[0]{\mathbb{P}_{k}(F)}^d\ \forall\ F \in \mathcal{F}_h,\ \bar{v}_h
      = 0 \text{ on } \Gamma_D},
    \\
    Q_h
    &:= \cbr[1]{q_h\in L^2(\Omega) : \ q_h \in \mathbb{P}_{k-1}(K) ,\
      \forall \ K \in \mathcal{T}_h},
    \\
    \bar{Q}_h
    &:= \cbr[1]{\bar{q}_h \in L^2(\Gamma_0) : \ \bar{q}_h \in
      \mathbb{P}_{k}(F) \ \forall\ F \in \mathcal{F}_h},
    \\
    \bar{Q}_h^0
    &:=\cbr[0]{\bar{\psi}_h \in \bar{Q}_h:
      \bar{\psi}_h=0 \text{ on } \Gamma_P}.    
  \end{split}
\end{equation}
Element and facet function pairs will be denoted by boldface, for
example,
\begin{equation*}
  \begin{split}
    \boldsymbol{v}_h &= (v_h, \bar{v}_h) \in \boldsymbol{V}_h := V_h\times \bar{V}_h,
    \qquad
    \boldsymbol{q}_h = (q_h, \bar{q}_h)\in \boldsymbol{Q}_h := Q_h\times \bar{Q}_h,
    \\
    \boldsymbol{\psi}_h &= (\psi_h, \bar{\psi}_h)\in \boldsymbol{Q}_h^0 := Q_h\times \bar{Q}_h^0,   
  \end{split}
\end{equation*}
and it will also be useful to define
$\boldsymbol{X}_h := \boldsymbol{V}_h \times \boldsymbol{Q}_h \times
V_h \times \boldsymbol{Q}_h^0$. 

\begin{remark}
  \label{rem:edghdg}
  The HDG method seeks an approximation in $\boldsymbol{X}_h$ with $V_h$, $\bar{V}_h$, $Q_h$, $\bar{Q}_h$, and $\bar{Q}_h^0$ defined in \cref{eq:approximationspaces}. If $\bar{V}_h$ is replaced by $\bar{V}_h \cap C^0(\Gamma_0)$ then we obtain the EDG-HDG method. The analysis in this paper holds for both the HDG and EDG-HDG methods. For notational purposes, in the analysis, $\boldsymbol{X}_h$ and $\boldsymbol{V}_h$ will refer both to the HDG and EDG-HDG spaces.
\end{remark}

For the analysis of the HDG and EDG-HDG methods we assume that the exact solution
is such that:
\begin{align*}
  u(t) &\in V:=\cbr[0]{v \in \sbr[0]{H^1(\Omega)}^d\,:\, v|_{\Gamma_D}=0} \cap \sbr[0]{H^2(\Omega)}^d,
  \\
  p_T(t) &\in Q:=H^1(\Omega),
  \\
  z(t) &\in Z:=\cbr[0]{v \in \sbr[0]{H^1(\Omega)}^d\,:\, v\cdot n = 0 \text{ on } \Gamma_F},
  \\
  p(t) &\in Q^0 :=\cbr[0]{q \in H^1(\Omega)\,:\, q|_{\Gamma_P} = 0 \text{ on } \Gamma_P} \cap H^2(\Omega).
\end{align*}
Denoting by $\bar{V}$, $\bar{Q}$, $\bar{Z}$, and $\bar{Q}^0$ the trace
spaces of, respectively, $V$, $Q$, $Z$, and $Q^0$ to the mesh
skeleton, we introduce the extended spaces
\begin{align*}
  \boldsymbol{V}(h)
  &:= \boldsymbol{V}_h + V \times \bar{V},
  &
    Z(h)
  &:= V_h + Z,
  \\
  \boldsymbol{Q}(h)
  &:= \boldsymbol{Q}_h + Q \times \bar{Q},
  &
    \boldsymbol{Q}^0(h)
  &:= \boldsymbol{Q}^0_h + Q^0 \times \bar{Q}^0.
\end{align*}
Norms on the extended spaces $\boldsymbol{V}(h)$, $\boldsymbol{Q}(h)$,
and $\boldsymbol{Q}^0(h)$ are defined as:
\begin{align*}
  \tnorm{\boldsymbol{v}}_{v}^2
  &:= \norm[0]{\varepsilon(v)}_{\Omega}^2 + \sum_{K\in \mathcal{T}_h} h_K^{-1}\norm[0]{v-\bar{v}}_{\partial K}^2
  & \forall \boldsymbol{v} \in \boldsymbol{V}(h),
  \\
  \tnorm{\boldsymbol{v}}_{v'}^2
  &:= \tnorm{\boldsymbol{v}}_v^2 + \sum_{K\in \mathcal{T}_h} h_K^{2}|v|_{2,K}^2,
  & \forall \boldsymbol{v} \in \boldsymbol{V}(h),
  \\
  \tnorm{\boldsymbol{q}_h}_{q}^2
  &:= \norm[0]{q_h}_{\Omega}^2 + \sum_{K\in\mathcal{T}_h} h_K \norm[0]{\bar{q}_h}_{\partial K}^2
  & \forall \boldsymbol{q}_h \in \boldsymbol{Q}(h).
\end{align*}

To conclude this section, we remark that $C > 0$ will denote a
constant independent of $h$ and the model parameters.

\subsection{The semi-discrete problem}
\label{sec:semidiscrete}

In this section, we present the semi-discrete problem and provide an
energy estimate for this discretization. The fully-discrete problem is
presented in \cref{ss:fullydiscrete} which is analysed in
\cref{sec:apriori}.

The semi-discrete HDG method for Biot's consolidation model
\cref{eq:biot,eq:bics} is given by: Find
$(\boldsymbol{u}_h, z_h) \in C^0(I; \boldsymbol{V}_h\times V_h)$ and
$(\boldsymbol{p}_{T,h}, \boldsymbol{p}_h) \in C^1(I; \boldsymbol{Q}_h
\times \boldsymbol{Q}_h^0)$ such that for all
$(\boldsymbol{v}_h,\boldsymbol{q}_{Th},w_h,\boldsymbol{q}_h) \in
\boldsymbol{V}_h \times \boldsymbol{Q}_h \times V_h \times
\boldsymbol{Q}_h^0$:
\begin{subequations}
  \label{eq:semi-discrete-hdg}
  \begin{align}
    \label{eq:semi-discrete-hdg_a}
    a_h(\boldsymbol{u}_h, \boldsymbol{v}_h)
    + b_h(\boldsymbol{v}_h,\boldsymbol{p}_{Th})
    &= \del[1]{f, v_h}_{\Omega},
    \\
    \label{eq:semi-discrete-hdg_b}
    b_h(\boldsymbol{u}_h, \boldsymbol{q}_{Th})
    - \del[1]{\lambda^{-1}(p_{Th}-\alpha p_h), q_{Th}}_{\Omega}
    &= 0,
    \\
    \label{eq:semi-discrete-hdg_c}
    \del[1]{\kappa^{-1}z_h, w_h}_{\Omega}
    + b_h((w_h,0), \boldsymbol{p}_h)
    &= 0,
    \\
    \label{eq:semi-discrete-hdg_d}
    \del[1]{ \partial_t \del[0]{c_0p_h + \lambda^{-1}\alpha\del[0]{\alpha p_h - p_{Th}}}, q_h}_{\Omega}
      - b_h((z_h,0), \boldsymbol{q}_h)
    &= \del[1]{g, q_h}_{\Omega},
  \end{align}
\end{subequations}
where
\begin{subequations}
  \begin{align}
    \label{eq:a_h}
    a_h(\boldsymbol{u}, \boldsymbol{v})
    :=&
        (2\mu \varepsilon(u), \varepsilon(v))_{\Omega}
        +\sum_{K\in\mathcal{T}}\langle\tfrac{2\beta\mu}{h_K}
        (u - \bar{u}),v-\bar{v}\rangle_{\partial K}
    \\ \nonumber
      &- \langle 2\mu \varepsilon(u)n,
        v-\bar{v}\rangle_{\partial \mathcal{T}}
        - \langle 2\mu \varepsilon(v)n,
        u-\bar{u}\rangle_{\partial \mathcal{T}},
    \\
    \label{eq:b_h}
    b_h(\boldsymbol{v}, \boldsymbol{q})
    :=&
        - (q, \nabla\cdot v )_{\Omega}
        + \langle \bar{q},
        (v-\bar{v}) \cdot n\rangle_{\partial \mathcal{T}}.
  \end{align}
\end{subequations}

To analyze the HDG and EDG-HDG methods, let us recall some properties of the
bilinear forms $a_h$ and $b_h$. It was shown in \cite[Lemma
4.2]{Rhebergen:2017} and \cite[Lemma 2]{Cesmelioglu:2020} that there
exist constants $C$ and $\beta_0 > 0$ such that for $\beta > \beta_0$,
\begin{equation}
  \label{eq:a_h-coercivity}
  a_h(\boldsymbol{v}_h, \boldsymbol{v}_h)
  \ge C \mu \tnorm{\boldsymbol{v}_h}_v^2 \qquad \forall \boldsymbol{v}_h \in \boldsymbol{V}_h.
\end{equation}
Additionally, $a_h$ satisfies the following continuity result
\cite[Lemma 3]{Cesmelioglu:2020}:
\begin{equation}
  \label{eq:a_h-continuity}
  a_h(\boldsymbol{u}, \boldsymbol{v}) \le C \mu \tnorm{\boldsymbol{u}}_{v'} \tnorm{\boldsymbol{v}}_{v'}
  \qquad \forall \boldsymbol{u}, \boldsymbol{v} \in \boldsymbol{V}(h).
\end{equation}
The bilinear form $b_h$ satisfies the following stability results:
\begin{subequations}
  \label{eq:bh-infsups}
  \begin{align}
    \label{eq:b_h-inf-sup}
    \inf_{\boldsymbol{q}_h \in \boldsymbol{Q}_h} \sup_{\boldsymbol{0}\neq\boldsymbol{v}_h \in \boldsymbol{V}_h}
    \frac{ b_h(\boldsymbol{v}_h, \boldsymbol{q}_h) }{\tnorm{\boldsymbol{v}_h}_{v}\tnorm{ \boldsymbol{q}_h }_{q}}
    & \ge C,
    \\
    \label{eq:b_h-inf-sup2}    
    \inf_{\boldsymbol{q}_h \in \boldsymbol{Q}_h^0} \sup_{0\neq w_h \in  V_h}
    \frac{b_h((w_h,0), \boldsymbol{q}_h)}{\norm[0]{w_h}_{\Omega}\tnorm{\boldsymbol{q}_h}_q}
    & \ge C,
  \end{align}
\end{subequations}
where the first inequality was shown in \cite[Lemma
1]{Rhebergen:2018b} and \cite[Lemma 8]{Rhebergen:2020} and the second is proven in
\cref{ap:inf-sup_bhwh0qh}. Continuity of the bilinear form $b_h$ was
established in \cite{Rhebergen:2017}:
\begin{equation}
  \label{eq:h_b-continuity}
  |b_h(\boldsymbol{v}, \boldsymbol{q})| \le C \tnorm{\boldsymbol{v}}_v \tnorm{\boldsymbol{q}}_q
  \qquad \forall \boldsymbol{v} \in \boldsymbol{V}(h), \boldsymbol{q} \in \boldsymbol{Q}(h).
\end{equation}

\begin{lemma}[Consistency]
  \label{lem:consistency}
  Let $(u, p_T, z, p)$ be a solution to \cref{eq:biot,eq:bics} and let
  $\bar{u}$, $\bar{p}_T$, and $\bar{p}$ be the traces of,
  respectively, $u$, $p_T$, and $p$ on the mesh skeleton. Then
  $(\boldsymbol{u}, \boldsymbol{p}_{T}, z, \boldsymbol{p})$ satisfies
  \cref{eq:semi-discrete-hdg}.
\end{lemma}
\begin{proof}
  The proof is standard and follows by integration by parts,
  smoothness of the solution to \cref{eq:biot,eq:bics},
  single-valuedness of $\bar{v}_h$ and $\bar{q}_h$ on element
  boundaries, and using that $\bar{v}_h=0$ on $\Gamma_D$ and
  $\bar{q}_h=0$ on $\Gamma_P$. 
\end{proof}
The following theorem now shows energy stability of the semi-discrete
problem.
\begin{theorem}[Stability]
  \label{thm:stability}
  Suppose that
  $(\boldsymbol{u}_h, \boldsymbol{p}_{Th}, z_h, \boldsymbol{p}_h) \in
  C^1(I; \boldsymbol{X}_h)$ is a solution to
  \cref{eq:semi-discrete-hdg} with $f \in W^{1,1}(I; L^2(\Omega))$
  and $g \in L^2(I; L^2(\Omega))$. Let $X(t) \ge 0$ and
  $Y(t) \ge 0$ be defined by:
  \begin{align*}
    X(t)^2
    =& a_h (\boldsymbol{u}_h(t), \boldsymbol{u}_h(t))
    \\
     &+ \del[1]{\lambda^{-1} \del[1]{p_{Th} (t) - \alpha p_h (t)}, p_{Th} (t) - \alpha p_h (t) }_{\Omega}
       + \del[1]{c_0 p_h (t), p_h (t)}_{\Omega},
    \\
    Y(t)^2
    =& \del[1]{\kappa^{-1} z_h(t), z_h(t) }_{\Omega}.
  \end{align*}
  Then, there exists $C>0$, independent of $t>0$, such that
  \begin{subequations}
    \label{eq:XY-estm}
    \begin{multline}
      \label{eq:XY-estm-a}
      X(t)
      \le X(0)
      + C \bigg[ \mu^{-1/2} \max_{0 \le s \le t} \norm[0]{f(s)}_{\Omega}
      \\
      + \mu^{-1/2} \int_0^t \norm[0]{\partial_tf(s)}_{\Omega} \dif s
      + \del[3]{\int_0^t \norm[0]{g(s)}_{\Omega}^2 \dif s}^{\frac{1}{2}} \bigg],      
    \end{multline}
    and
    \begin{multline}
      \label{eq:XY-estm-b}
      \del[3]{ \int_0^t Y(s)^2 \dif s }^{\frac{1}{2}}
      \le C \bigg[ X(0)
      + \mu^{-1/2}\max_{0 \le s \le t} \norm[0]{f(s)}_{\Omega}
      \\
      + \mu^{-1/2}\int_0^t \norm[0]{\partial_tf(s)}_{\Omega} \dif s
      + \del[3]{ \int_0^t \norm[0]{g(s)}_{\Omega}^2 \dif s }^{\frac{1}{2}} \bigg].                  
    \end{multline}
  \end{subequations}
\end{theorem}
\begin{proof}
  We first note that by the inf-sup condition \cref{eq:b_h-inf-sup2},
  \cref{eq:semi-discrete-hdg_c}, and the Cauchy--Schwarz inequality,
  \begin{equation}
    \label{eq:p-z-estm}
    \begin{split}
      \norm{p_h(t)}_{\Omega}
      &\le C \sup_{0\neq w_h \in  V_h}
      \frac{|b_h((w_h,0), \boldsymbol{p}_h)|}{\norm[0]{w_h}_{\Omega}}
      \\
      &\le C \sup_{0\neq w_h \in  V_h}
      \frac{|(\kappa^{-1}z_h, w_h)|}{\norm[0]{w_h}_{\Omega}}
      \le C \kappa^{-1}\norm[0]{z_h(t)}_{\Omega}.      
    \end{split}
  \end{equation}
  Now, in
  \cref{eq:semi-discrete-hdg_a,eq:semi-discrete-hdg_c,eq:semi-discrete-hdg_d}
  set
  $(\boldsymbol{v}_h, w_h, \boldsymbol{q}_h) =
  (\partial_t\boldsymbol{u}_h, z_h, \boldsymbol{p}_h)$. Take the time
  derivative of \cref{eq:semi-discrete-hdg_b} and set
  $\boldsymbol{q}_{Th} = -\boldsymbol{p}_{Th}$. Adding the resulting
  equations we find:
  \begin{equation}
    \label{eq:energy-eq}
    \frac{1}{2} \od{}{t} X(t)^2 + Y(t)^2 = \del[0]{f(t), \partial_tu_h(t)}_{\Omega} + (g(t), p_h(t))_{\Omega}.
  \end{equation}
  Integrating \cref{eq:energy-eq} in time from 0 to $t$ results in
  \begin{equation*}
    \frac{1}{2} ( X(t)^2 - X(0)^2) + \int_0^t Y(s)^2 \dif s
    = \int_0^t \sbr[1]{ \del[0]{f(s), \partial_t u_h(s)}_{\Omega} + (g(s), p_h(s))_{\Omega} } \dif s.
  \end{equation*}
  Integration by parts, Young's inequality and \cref{eq:p-z-estm}
  imply
  \begin{equation*}
    \begin{split}
      \frac{1}{2} ( X(t)^2 - X(0)^2) + \int_0^t Y(s)^2 \dif s
      \le& \del[0]{f(t), u_h(t)}_{\Omega} - \del[0]{f(0), u_h(0)}_{\Omega}
      \\
      & - \int_0^t \del[0]{\partial_tf(s), u_h(s)}_{\Omega} \dif s
      \\
      &+ C \int_0^t\norm[0]{g(s)}_{\Omega}^2 \dif s + \frac{1}{2} \int_0^t Y(s)^2 \dif s.      
    \end{split}
  \end{equation*}
  Coercivity of $a_h$ \cref{eq:a_h-coercivity} and a discrete Korn's
  inequality imply that
  \begin{equation*}
    \norm[0]{u_h(t)}_{\Omega} \le C\mu^{-1/2} a_h(\boldsymbol{u}_h, \boldsymbol{u}_h)^{1/2} \le C\mu^{-1/2} X(t).
  \end{equation*}
  Therefore, by the Cauchy-Schwarz inequality, for any $t \ge 0$,
  \begin{equation}
    \begin{split}
      \label{eq:XY-ineq}  
      X(t)^2 + \int_0^t Y(s)^2 \dif s
      \le& X(0)^2 + C\mu^{-1/2} \del[2]{ \norm{f(t)}_{\Omega} X(t) + \norm{f(0)}_{\Omega} X(0) }
      \\
      &+ C\mu^{-1/2} \int_0^t \norm[0]{\partial_t f(s)}_{\Omega} X(s) \dif s + C \int_0^t \norm{g(s)}_{\Omega}^2 \dif s.      
    \end{split}
  \end{equation}
  To obtain \cref{eq:XY-estm-a}, we may assume, without loss of generality, that
  \begin{align}
    \label{eq:X-max-assumption}
    \max_{0 \le s \le t} X(s) = X(t)>0. 
  \end{align}
  Note that if \cref{eq:X-max-assumption} does not hold, then there
  exists a $t_M$ such that $\max_{0 \le s \le t} X(s) = X(t_M)$ for
  $0 \le t_M < t$. The estimate \cref{eq:XY-estm-a} for $X(t_M)$ then
  implies \cref{eq:XY-estm-a} for $X(t)$. From
  \cref{eq:XY-ineq,eq:X-max-assumption}, we then find
  \begin{equation}
    \label{eq:stability-1}
    \begin{split}
      &X(t)^2 + \int_0^t Y(s)^2 \dif s \\
      \le & \del[2]{X(0) + C \mu^{-1/2} \del[1]{\, \norm[0]{f(t)}_{\Omega}  + \norm[0]{f(0)}_{\Omega}
          + \int_0^t \norm[0]{\partial_t f(s)}_{\Omega} \dif s}} X(t) 
      \\
      & + C \int_0^t  \norm[0]{g(s)}_{\Omega}^2 \dif s.
    \end{split}    
  \end{equation}
  Define
  $\alpha(t):=\del[1]{C\int_0^t \norm[0]{g(s)}_{\Omega}^2
    \dif s}^{1/2} > 0$. If $\alpha (t) \leq X(t)$, dividing
  \cref{eq:stability-1} by $X(t)$ implies
  \begin{equation*}
    X(t) \le X(0) + C \mu^{-1/2} \del[1]{\, \norm[0]{f(t)}_{\Omega}  + \norm[0]{f(0)}_{\Omega}
      + \int_0^t \norm[0]{\partial_tf(s)}_{\Omega} \dif s} + \alpha(t).
  \end{equation*}
  Note that this inequality holds trivially if $X(t) <
  \alpha(t)$. Proceeding, we find
  {\fontsize{9.5}{10.5}\selectfont
  \begin{align*}
    X(t)
    &\le X(0) + C\mu^{-1/2} \del[2]{ 2\max_{0\leq s\leq t}\norm[0]{f(s)}_{\Omega}  + \int_0^t \norm[0]{\partial_tf(s)}_{\Omega} \dif s} + \alpha(t)
    \\
     &= X(0) + C\mu^{-1/2} \del[2]{ 2\max_{0\leq s\leq t}\norm[0]{f(s)}_{\Omega}  + \int_0^t \norm[0]{\partial_tf(s)}_{\Omega} \dif s}
       + C \del[2]{\int_0^t \norm[0]{g(s)}_{\Omega}^2 \dif s }^{1/2},
  \end{align*}
  }
  so that \cref{eq:XY-estm-a} follows. \Cref{eq:XY-estm-b} follows by
  combining \cref{eq:XY-estm-a} and \cref{eq:stability-1}. 
\end{proof}

\subsection{The fully discrete problem}
\label{ss:fullydiscrete}
To define the fully discrete scheme, let
$\cbr[0]{t^n}_{0\leq n\leq N}$ be a uniform partition of $I$ and let
$\Delta t > 0$ be the corresponding time step. We will denote the
value of a function $f(t)$ at $t=t^n$ by $f^n := f(t^n)$. For a
sequence $\cbr[0]{f^n}_{n\ge 1}$,
$d_tf^{n} := \frac{f^{n} - f^{n-1}}{\Delta t}$ defines a first order
difference operator. Note that we use the superscript $n$ to denote the time level. This is not to be confused with the normal vector $n$. Using Backward Euler time stepping, the fully
discrete problem reads: Find
$(\boldsymbol{u}_h^{n+1},\boldsymbol{p}_{Th}^{n+1}, z_h^{n+1},
\boldsymbol{p}_h^{n+1}) \in \boldsymbol{X}_h$, with $n\geq 0$, such
that
\begin{subequations}
  \label{eq:fully-discrete-hdg}
  \begin{align}
    \label{eq:fully-discrete-hdg_a}
    & a_h(\boldsymbol{u}_h^{n+1}, \boldsymbol{v}_h)
    + b_h(\boldsymbol{v}_h,\boldsymbol{p}_{Th}^{n+1})
    = \del[1]{f^{n+1}, v_h}_{\Omega},
    \\
    \label{eq:fully-discrete-hdg_b}
    & b_h(\boldsymbol{u}_h^{n+1}, \boldsymbol{q}_{Th})
    + \lambda^{-1}\del[1]{\alpha p_{h}^{n+1}-p_{Th}^{n+1}, q_{Th}}_{\Omega}
    = 0,
    \\
    \label{eq:fully-discrete-hdg_c}
    & \del[1]{\kappa^{-1}z_h^{n+1}, w_h}_{\Omega}
    + b_h((w_h,0), \boldsymbol{p}_h^{n+1})
    = 0,
    \\
    \label{eq:fully-discrete-hdg_d}
    &\tfrac1{\Delta t}\del[1]{(c_0 p_h^{n+1}, q_h)_{\Omega} + \lambda^{-1}(\alpha p_{h}^{n+1}-p_{Th}^{n+1},\alpha q_h)_{\Omega}}
    - b_h((z_h^{n+1},0), \boldsymbol{q}_h),
    \\
    \nonumber
    &= \tfrac1{\Delta t}\del[1]{(c_0p_h^{n}, q_h)_{\Omega}+ \lambda^{-1}(\alpha p_h^{n}-p_{Th}^{n},\alpha q_h)_{\Omega}}
    + \del[1]{g^{n+1}, q_h}_{\Omega},
  \end{align}
\end{subequations}
for all
$(\boldsymbol{v}_h,\boldsymbol{q}_{Th}, w_h, \boldsymbol{q}_h) \in
\boldsymbol{X}_h$. We first show that
\cref{eq:fully-discrete-hdg} is well-posed.

\begin{theorem}
  There exists a unique solution to \cref{eq:fully-discrete-hdg}.
\end{theorem}
\begin{proof}
  It is sufficient to show that if the data is equal to zero then the solution is
  zero. As such, suppose that $f^{n+1}=0$, $g^{n+1}=0$,
  $\boldsymbol{p}_{Th}^n=\boldsymbol{0}$, and
  $\boldsymbol{p}_h^n=\boldsymbol{0}$. Then, setting
  $\boldsymbol{v}_h=\boldsymbol{u}_h^{n+1}$,
  $\boldsymbol{q}_{Th}=-\boldsymbol{p}_{Th}^{n+1}$,
  $w_h=\Delta tz_h^{n+1}$, and
  $\boldsymbol{q}_{h}=\Delta t\boldsymbol{p}_{h}^{n+1}$ in
  \cref{eq:fully-discrete-hdg} and adding the equations, we obtain:
  \begin{equation*}
    a_h(\boldsymbol{u}_h^{n+1}, \boldsymbol{u}_h^{n+1})
    + c_0\norm[0]{p_h^{n+1}}_{\Omega}^2
    +\lambda^{-1}\norm[0]{\alpha p_{h}^{n+1}-p_{Th}^{n+1}}_{\Omega}^2
    +\kappa^{-1} \Delta t \norm[0]{z_h^{n+1}}_{\Omega}^2
    =0.
  \end{equation*}
  Coercivity of $a_h$ \cref{eq:a_h-coercivity}, positivity of $\kappa$
  and $\lambda$, and nonnegativity of $c_0$ directly imply that
  $\boldsymbol{u}_h^{n+1}=\boldsymbol{0}$ and
  $z_h^{n+1}=0$. Substituting $\boldsymbol{u}_h^{n+1}=\boldsymbol{0}$
  in \cref{eq:fully-discrete-hdg_a}, $\boldsymbol{p}_{Th}^{n+1}=0$
  follows from the inf-sup condition \cref{eq:b_h-inf-sup}. This then
  implies $p_{h}^{n+1}=0$ since $\alpha, \lambda>0$.

  Using a BDM local lifting of the normal trace \cite[Proposition
  2.10]{Du:book}, there exists $\tilde{w}_h\in V_h$ such that
  $\langle\tilde{w}_h\cdot n,
  \bar{p}_h^{n+1}\rangle_{\partial\mathcal{T}_h}=
  \norm[0]{\bar{p}_h^{n+1}}^2_{\partial\mathcal{T}_h}$. Setting
  $z_h^{n+1}=0$, $p_h^{n+1}=0$ and choosing $w_h=\tilde{w}_h$ in
  \cref{eq:fully-discrete-hdg_c}, we obtain $\bar{p}_h^{n+1}=0$. This
  completes the proof. 
\end{proof}

Let us also note that the fully-discrete scheme
\cref{eq:fully-discrete-hdg} results in divergence-conforming
solutions for the displacement $u_h^n$ and velocity $z_h^n$. To see
this, set $\boldsymbol{v}_h = \boldsymbol{0}$,
$\boldsymbol{q}_h = \boldsymbol{0}$, $w_h=0$, and $q_{Th} = 0$ in
\cref{eq:semi-discrete-hdg} and note that since
$u_h^n \cdot n \in P_k(F)$ and $\bar{u}_h^n=0$ on $\Gamma_D$,
\begin{subequations}
  \label{eq:div-conf-u}
  \begin{align}
    \jump{u_h^n\cdot n} &= 0, &&\forall x\in F,\quad \forall F\in \mathcal{F}\backslash\mathcal{F}_T,
    \\
    u_h^n\cdot n &= \bar{u}_h^n\cdot n, &&\forall x\in F,\quad \forall F\in \mathcal{F}_T.
  \end{align}  
\end{subequations}
Similarly, by setting $\boldsymbol{v}_h = \boldsymbol{0}$,
$\boldsymbol{q}_{Th} = \boldsymbol{0}$, $w_h=0$, and $q_h = 0$ in
\cref{eq:semi-discrete-hdg}, and noting that
$z_h^n \cdot n \in P_k(F)$, we find that
\begin{equation}
  \label{eq:div-conf-z}
  \jump{z_h^n\cdot n} = 0, \quad \forall x\in F,\quad \forall F\in \mathcal{F}\backslash \mathcal{F}_P.
\end{equation}

\section{A priori error estimates}
\label{sec:apriori}
To facilitate the a priori error analysis, we introduce various
interpolation operators. First, let
$\Pi_V: \sbr[0]{H^1(\Omega)}^d \rightarrow V_h$ be the BDM
interpolation operator \cite[Section III.3]{Brezzi:book}, \cite[Lemma
7]{Hansbo:2002} with the following interpolation estimate:
\begin{equation}
  \label{eq:interpolation-property}
  \norm[0]{z - \Pi_V z}_{K} \leq Ch_K^{\ell} \norm[0]{z}_{\ell,K},
  \quad 1 \leq \ell \leq k+1.
\end{equation}
The elliptic interpolation operator,
$\boldsymbol{\Pi}_{V}^{\text{ell}} := (\Pi_V^{\text{ell}},
\bar{\Pi}_V^{\text{ell}}): \sbr[0]{H^1(\Omega)}^d\rightarrow
\boldsymbol{V}_h$ is defined by:
\begin{equation*}
  a_h(\boldsymbol{\Pi}_{V}^{\text{ell}}u, \boldsymbol{v}_h) 
  =a_h((u,u), \boldsymbol{v}_h), \qquad \forall \boldsymbol{v}_h \in \boldsymbol{V}_h.
\end{equation*}
Standard a priori error estimate theory for second order elliptic
equations imply
\begin{equation}
  \label{eq:ell-interpolation-property}
  a_h( \boldsymbol{u} - \boldsymbol{\Pi}_V^{\text{ell}} u, \boldsymbol{u} - \boldsymbol{\Pi}_V^{\text{ell}} u )^{\frac 12}
  \leq C \mu^{1/2} h_K^{\ell-1}\norm[0]{u}_{\ell,\Omega} \quad 1\leq \ell\leq k+1.
\end{equation}
By $\Pi_Q$, $\bar{\Pi}_{Q}$, and $\bar{\Pi}_{Q^0}$ we
denote the $L^2$-projections onto, respectively, $Q_h$ and the trace
spaces $\bar{Q}_h$ and $\bar{Q}_h^0$. Given the
interpolation/projection operators, the numerical initial data is set
by imposing the interpolation/projection of continuous initial data as
follows:
\begin{multline}
  \label{eq:initial-data}
  (({u}_h^0, \bar{u}_h^0), (p_{Th}^0, \bar{p}_{Th}^0), z_h^0, (p_h^0, \bar{p}_h^0)) \\
  = ((\Pi_V^{\text{ell}} u(0), \bar{\Pi}_V^{\text{ell}} u(0)), (\Pi_Q p_{Th} (0), \bar{\Pi}_Q p_{Th}(0)), \Pi_V z(0), (\Pi_{Q^0} p(0), \bar{\Pi}_{Q^0} p(0))) .
\end{multline}
%
\begin{equation}
  \label{eq:initial-data}
  \boldsymbol{u}_h^0=({u}_h^0, \bar{u}_h^0) = (\Pi_V^{\text{ell}} u(0), \bar{\Pi}_V^{\text{ell}} u(0)),\quad  p_h^0= \Pi_{Q^0} p(0).
\end{equation}

In the error analysis it will be convenient to split the error into
approximation and interpolation errors:
\begin{subequations}
  \label{eq:errorsplit}
  \begin{align}
    \omega- \omega_h
    &= e_{\omega}^I - e_{\omega}^h, & \omega &= u, p_T, z, p,
    \\
    \xi|_{\Gamma^0} - \bar{\xi}_h
    &= \bar{e}_{\xi}^I - \bar{e}_{\xi}^h, & \xi &= u, p_T, p,  
  \end{align}  
\end{subequations}
where
\begin{align*}
  e_{u}^I &= u - \Pi_{V}^{\text{ell}} u, & e_{u}^h &= u_h - \Pi_{V}^{\text{ell}} u, & e_{z}^I &= z - \Pi_{V}z, & e_{z}^h &= z_h - \Pi_{V}z,
  \\
  e_p^I &= p - \Pi_Q p, & e_p^h &= p_h - \Pi_Q p, & e_{p_T}^I &= p_T - \Pi_Q p_T, & e_{p_T}^h &= p_{Th} - \Pi_Q p_T,
\end{align*}
and where
\begin{align*}
  \bar{e}_{u}^I  &= u|_{\Gamma^0} - \bar{\Pi}_V^{\text{ell}} u,
  &\bar{e}_{u}^h  &= \bar{u}_h - \bar{\Pi}_V^{\text{ell}} u,
  &\bar{e}_{p}^I &= p|_{\Gamma^0} - \bar{\Pi}_{Q^0} p,
  & \bar{e}_{p}^h &= \bar{p}_h - \bar{\Pi}_{Q^0}p, 
  \\
  \bar{e}_{p_T}^I  &= p_T|_{\Gamma^0} - \bar{\Pi}_Q p_T, & \bar{e}_{p_T}^h &= \bar{p}_{Th} - \bar{\Pi}_Q p_T. &&  
\end{align*}
Following the convention introduced earlier in this paper, we use
boldface notation for element/facet error pairs, i.e.,
$\boldsymbol{e}_{\xi}^I = (e_{\xi}^I, \bar{e}_{\xi}^I)$ and
$\boldsymbol{e}_{\xi}^h = (e_{\xi}^h, \bar{e}_{\xi}^h)$ for
$\xi = u, p_T, p$.

It will also be useful to introduce the following error estimates: let
$\psi$ be a regular enough function defined on $[0,T] \times D$, for
some domain $D \subset \mathbb{R}^d$, then, as a consequence of Taylor's
theorem (see \cref{ap:err-taylor}),
\begin{subequations}
  \begin{align}
    \label{ineq:taylor3}
    \sum_{m=1}^n \Delta t  \|\partial_t \psi^m-d_t\psi^m\|_{0,D}
    &\le \Delta t \| \partial_{tt} \psi \|_{L^1(I; L^2(D))}, &&
    \\
    \label{ineq:taylor4}
    \Delta t \sum_{m=1}^n \| d_t e_{\psi}^{I,m} \|_{0,D} 
    &\le C h^l \| \partial_t \psi \|_{L^1(I; H^l(D))},
    && \psi = p,p_T.
  \end{align}
\end{subequations}

\begin{lemma}
  \label{lemma:not-gronwall}
  Let $\{A_i\}_i, \{B_i\}_i$, $\{E_i\}_i$, and $\{D_i\}_i$ be
  nonnegative sequences. Suppose these sequences satisfy
  \begin{equation}
    \label{eq:abcd-assumption}
    A_n^2 + \sum_{i=0}^n B_i^2 \leq  A_0^2 + \sum_{i=1}^{n} E_i A_i + \sum_{i=0}^n D_i,
  \end{equation}	
  for all $n \geq 0$. Then for any $n\geq 0$,
  \begin{subequations}
    \begin{align}
      \label{eq:an-estimate}
      A_n
      &\leq  A_0 +  \sum_{i=1}^{n} E_i + \del[2]{\sum_{i=0}^n D_i}^{1/2} ,
      \\
      \label{eq:b_sum-estimate}
      \del[2]{\sum_{i=0}^n B_i^2}^{1/2}
      &\leq C \del[3]{A_0 +  \sum_{i=1}^{n} E_i + \del[2]{\sum_{i=0}^n D_i }^{1/2}},
    \end{align}	    
  \end{subequations}
  with $C>0$ independent of $n$.
\end{lemma}
\begin{proof}
  First, note that \cref{eq:abcd-assumption} and \cref{eq:an-estimate}
  directly imply \cref{eq:b_sum-estimate}. It is therefore sufficient
  to prove \cref{eq:an-estimate}. Similar to the assumption made in
  the proof of \cref{thm:stability}, we assume without loss of
  generality that $0 < A_n = \max_{0 \le i \le n} A_i$.  If
  $A_n\leq \del[1]{\sum_{i=0}^n D_i}^{1/2}$, then
  \cref{eq:an-estimate} is satisfied trivially. On the other hand, if
  $A_n> \del[1]{\sum_{i=0}^n D_i}^{1/2}$, then
  \cref{eq:abcd-assumption} implies
  \begin{equation*}
    A_n^2 + \sum_{i=0}^n B_i^2 \leq  A_0A_n + A_n \sum_{i=0}^{n} E_i  + A_n\Big(\sum_{i=0}^n D_i\Big)^{1/2}.
  \end{equation*}
  The result now follows by dividing by $A_n$. 
\end{proof}

Before addressing the main result (\cref{thm:error-estimate}) of this
section, we first determine the error equation.

\begin{lemma}[Error equation]
  \label{lem:errorequation} Suppose that $\{(\boldsymbol{u}_h^n, \boldsymbol{p}_{Th}^n, z_h^n, \boldsymbol{p}_h^n)\}_{n \ge 1}$ is the solution of \eqref{eq:fully-discrete-hdg} with the numerical initial data \eqref{eq:initial-data}.
  The approximation and interpolation errors satisfy
  \begin{align*}
    &a_h(\boldsymbol{e}_u^{h,n+1},\boldsymbol{v}_h) + b_h(\boldsymbol{v}_h, \boldsymbol{e}_{p_{T}}^{h,n+1})
    \\
    &
      + b_h(d_t\boldsymbol{e}_u^{h,n+1}, \boldsymbol{q}_{Th}) + \lambda^{-1}(d_t(\alpha e_p^{h,n+1}-e_{p_T}^{h,n+1}), \alpha q_h+q_{Th})_{\Omega} \\
    &+ \kappa^{-1}(e_z^{h,n+1}, w_h)_{\Omega} - b_h((e_z^{h,n+1},0),\boldsymbol{q}_h)
    \\
    &+ b_h((w_h,0),\boldsymbol{e}_p^{h,n+1})_{\Omega} + (c_0 d_t e_p^{h,n+1}, q_h)_{\Omega} 
    \\
    =& b_h(d_t\boldsymbol{e}_u^{I,n+1}, \boldsymbol{q}_{Th}) + \kappa^{-1}(e_z^{I,n+1}, w_h)_{\Omega} + (c_0(\partial_tp^{n+1}-d_t \Pi_Q p^{n+1}), q_h)_{\Omega}  \\
    &+ \lambda^{-1}(\alpha (\partial_tp^{n+1} - d_t \Pi_Q p^{n+1}) -(\partial_tp_T^{n+1}-d_t\Pi_Q p_T^{n+1}), \alpha q_h)_{\Omega},
  \end{align*}
  for any
  $(\boldsymbol{v}_h,\boldsymbol{q}_{Th}, w_h, \boldsymbol{q}_h) \in
  \boldsymbol{X}_h$.
\end{lemma}
\begin{proof}  
  By \cref{lem:consistency}, we can substitute $(u^{n+1}, p_T^{n+1}, z^{n+1}, p^{n+1})$, the solution of
  \cref{eq:biot,eq:bics} evaluated at $t=t^{n+1}$, into
  \cref{eq:semi-discrete-hdg}. Then, 
  subtracting \cref{eq:fully-discrete-hdg}, applying $d_t$ to the
  second equation of the resulting set of equations, and adding all
  the equations, we obtain using \cref{eq:errorsplit}:
  \begin{align*}
    &a_h(\boldsymbol{e}_u^{h,n+1},\boldsymbol{v}_h) + b_h(\boldsymbol{v}_h, \boldsymbol{e}_{p_{T}}^{h,n+1})
    \\
    &
      + b_h(d_t\boldsymbol{e}_u^{h,n+1}, \boldsymbol{q}_{Th}) + \lambda^{-1}(d_t(\alpha e_p^{h,n+1}-e_{p_T}^{h,n+1}), \alpha q_h + q_{Th})_{\Omega} \\
    &+ \kappa^{-1}(e_z^{h,n+1}, w_h)_{\Omega} - b_h((e_z^{h,n+1},0),\boldsymbol{q}_h)
    \\
    &+ b_h((w_h,0),\boldsymbol{e}_p^{h,n+1})_{\Omega} + (c_0 d_t e_p^{h,n+1}, q_h )_{\Omega} 
    \\
    =& a_h(\boldsymbol{e}_u^{I,n+1},\boldsymbol{v}_h) + b_h(\boldsymbol{v}_h, \boldsymbol{e}_{p_T}^{I,n+1}) \\
    &+ b_h(d_t\boldsymbol{e}_u^{I,n+1}, \boldsymbol{q}_{Th}) + \lambda^{-1}(d_t(\alpha e_p^{I,n+1}-e_{p_T}^{I,n+1}), \alpha q_h + q_{Th})_{\Omega} \\
    &+ \kappa^{-1}(e_z^{I,n+1}, w_h)_{\Omega} - b_h((e_z^{I,n+1},0),\boldsymbol{q}_h) \\
    &+ b_h((w_h,0),\boldsymbol{e}_p^{I,n+1}) + (c_0(\partial_tp^{n+1}-d_t \Pi_Q p^{n+1}), q_h)_{\Omega}  \\
    &+ \lambda^{-1}(\alpha (\partial_tp^{n+1} - d_t \Pi_Q p^{n+1}) -(\partial_tp_T^{n+1}-d_t\Pi_Q p_T^{n+1}), \alpha q_h )_{\Omega}.
  \end{align*}
  The result follows by noting that:
  $a_h(\boldsymbol{e}_u^{I,n+1},\boldsymbol{v}_h) = 0$ by definition
  of $\boldsymbol{\Pi}_V^{\text{ell}}$;
  $b_h(\boldsymbol{v}_h, \boldsymbol{e}_{p_T}^{I,n+1}) = b_h((w_h,0),\boldsymbol{e}_p^{I,n+1}) = 0$ because
  $\Pi_Q$ and $\bar{\Pi}_Q$ are $L^2$ projections into $Q_h$ and
  $\bar{Q}_h$, respectively, and $\nabla \cdot V_h = Q_h$;
  $b_h((e_z^{I,n+1},0),\boldsymbol{q}_h) = 0$ by the commuting
  property of the BDM interpolation operator
  $\Pi_Q \nabla \cdot {v} = \nabla \cdot \Pi_V {v}$ for
  $v \in \sbr[0]{H^1(\Omega)}^d$, the $H(\text{div})$-conformity of $e_z^{I,n+1}$, and the boundary conditions on $\Gamma_F$; and
  $\lambda^{-1}(d_t(\alpha e_p^{I,n+1}-e_{p_T}^{I,n+1}), \alpha q_h +
  q_{Th})_{\Omega} = 0$ because $\Pi_Q$ is the $L^2$-projection into
  $Q_h$. 
\end{proof}

We are now ready to prove an a priori error estimate for the HDG and EDG-HDG
methods \cref{eq:fully-discrete-hdg}. 

\begin{theorem}
  \label{thm:error-estimate}
  Let $(u, p_T, z, p)$ be a solution to \cref{eq:biot,eq:bics} on the time interval $I=(0,T]$ and let
  $\bar{u}$, $\bar{p}_T$, and $\bar{p}$ be the traces of,
  respectively, $u$, $p_T$, and $p$ on the mesh skeleton. Let
  $(\boldsymbol{u}_h^{n},\boldsymbol{p}_{Th}^{n}, z_h^{n},
  \boldsymbol{p}_h^{n}) \in \boldsymbol{X}_h$ be the solution to
  \cref{eq:fully-discrete-hdg}. Suppose the numerical initial data is
  imposed according to \cref{eq:initial-data}. The following error
  estimates hold: 
  \begin{subequations}
    \begin{align}
      c_0^{1/2}\norm[0]{p^n-p_h^n}_{\Omega} + \lambda^{-1/2}\norm[0]{\alpha (p^n-p_h^n) - (p_T^n-p_{Th}^n)}_{\Omega} &
      \nonumber
        \\
        + \mu^{1/2}\tnorm{\boldsymbol{u}^n-\boldsymbol{u}_h^{n}}_v
        + \kappa^{-1/2} \del[2]{\sum_{i=1}^n\Delta t\norm[0]{z^i-z_h^i}_{\Omega}^2}^{1/2}
      &\le C_1 \Delta t + C_2 h^l,      \label{eq:error-estimate1-a}
      \\
      \label{eq:error-estimate1-b}
      \mu^{-1/2}\tnorm{\boldsymbol{p}_T^n - \boldsymbol{p}_{Th}^n}_{q}
      &\le C_1 \Delta t + C_3 h^l,
    \end{align}
  \end{subequations}
  where
  \begin{align*}
    C_1
    =& C\max\cbr[0]{c_0^{1/2}, \lambda^{-1/2}}\norm[0]{p,p_T}_{W^{2,1}(I;L^2(\Omega))},
    \\
    C_2
    =& C \big(\mu^{1/2} \norm[0]{u}_{W^{1,1}(I; H^{l+1}(\Omega))}
    + \kappa^{-1/2} T^{1/2} \| z \|_{C^0(I; H^l(\Omega))}
    \\
     &+ \max\cbr[0]{c_0^{1/2}, \lambda^{-1/2}}\norm[0]{p,p_T}_{W^{1,1}(I; H^l(\Omega))}\big),
    \\
    C_3
    =&    C \big(\mu^{1/2} \norm{\partial_tu}_{L^1(I;H^{l+1}(\Omega))}
    + \max\cbr[0]{c_0^{1/2}, \mu^{-1/2}}\norm[0]{p,p_T}_{W^{1,1}(I;H^l(\Omega))}
    \\
    &+ \kappa^{-1/2} T^{1/2} \| z \|_{C^0(I; H^l(\Omega))}\big).
  \end{align*}
\end{theorem}
\begin{proof}
  Choose $\boldsymbol{v}_h = d_t\boldsymbol{e}_u^{h,n+1}$,
  $\boldsymbol{q}_{Th}=-\boldsymbol{e}_{p_T}^{h,n+1}$,
  $w_h=e_z^{h,n+1}$, and $\boldsymbol{q}_h = \boldsymbol{e}_p^{h,n+1}$
  in the error equation in \cref{lem:errorequation}. Then,
  \begin{align*}
    &a_h(\boldsymbol{e}_u^{h,n+1},d_t\boldsymbol{e}_u^{h,n+1}) + \kappa^{-1}(e_z^{h,n+1}, e_z^{h,n+1})_{\Omega}
    \\
    &+(c_0 d_t e_p^{h,n+1}, e_p^{h,n+1})_{\Omega} + \lambda^{-1}(d_t(\alpha e_p^{h,n+1}-e_{p_T}^{h,n+1}), \alpha e_p^{h,n+1}-e_{p_T}^{h,n+1})_{\Omega}
    \\
    = & - b_h(d_t\boldsymbol{e}_u^{I,n+1},\boldsymbol{e}_{p_T}^{h,n+1}) + \kappa^{-1}(e_z^{I,n+1}, e_z^{h,n+1})_{\Omega}
    \\
    &+ (c_0(\partial_tp^{n+1}-d_t \Pi_Q p^{n+1}), e_p^{h,n+1})_{\Omega}
    \\
    & + \lambda^{-1}(\alpha (\partial_tp^{n+1} - d_t \Pi_Q p^{n+1}) -(\partial_tp_T^{n+1}-d_t\Pi_Q p_T^{n+1}), \alpha e_p^{h,n+1})_{\Omega}.
\end{align*}
Using
$a(a-b) = \tfrac{1}{2}(a^2 + (b-a)^2 - b^2) \ge \tfrac{1}{2}(a^2 -
b^2)$ and multiplying both sides of the resulting inequality by
$\Delta t$, we arrive at
\begin{equation}
  \label{eq:error-eqs-sum}
  \begin{split}
  & \frac{c_0}2 \del[2]{\|e_p^{h,n+1}\|_{\Omega}^2 - \|e_p^{h,n}\|_{\Omega}^2}
  + \frac{\lambda^{-1}}2\del[2]{\|\alpha e_p^{h,n+1}-e_{p_T}^{h,n+1}\|_{\Omega}^2-\|\alpha e_p^{h,n}-e_{p_T}^{h,n}\|_{\Omega}^2}
  \\
  &+\frac 12 ( a_h(\boldsymbol{e}_u^{h,n+1},\boldsymbol{e}_u^{h,n+1}) -  a_h(\boldsymbol{e}_u^{h,n}, \boldsymbol{e}_u^{h,n}) )
  + \kappa^{-1}\Delta t \|e_z^{h,n+1}\|_{\Omega}^2
  \\
  \le& - \Delta t b_h(d_t\boldsymbol{e}_u^{I,n+1},\boldsymbol{e}_{p_T}^{h,n+1}) + \Delta t \kappa^{-1} (e_z^{I,n+1}, e_z^{h,n+1})_{\Omega}
  \\
  & + \Delta t (c_0(\partial_tp^{n+1}-d_t \Pi_Q p^{n+1}), e_p^{h,n+1})_{\Omega}
  \\
  & + \lambda^{-1}\Delta t (\alpha (\partial_tp^{n+1} - d_t \Pi_Q p^{n+1}) -(\partial_tp_T^{n+1}-d_t\Pi_Q p_T^{n+1}), \alpha e_p^{h,n+1})_{\Omega}
  \\
  =:& I_1^{n+1} + I_2^{n+1} + I_3^{n+1} + I_4^{n+1}.
  \end{split}
\end{equation}
We define
\begin{equation*}
  \begin{split}
    A_i^2
    &:= \frac{c_0}2 \|e_p^{h,i}\|_{\Omega}^2 
    + \frac{\lambda^{-1}}2\|\alpha e_p^{h,i}-e_{p_T}^{h,i}\|_{\Omega}^2 + \frac 12 a_h (\boldsymbol{e}_u^{h,i},\boldsymbol{e}_u^{h,i}),
    \\
    B_i^2
    &:= \frac{\kappa^{-1}}2 \Delta t \|e_z^{h,i}\|_{\Omega}^2,    
  \end{split}
\end{equation*}
so that \cref{eq:error-eqs-sum} can be written as
\begin{equation}
  \label{eq:An-Bn-ineq0}
  A_{n+1}^2 + 2 B_{n+1}^2 \le A_n^2 + I_1^{n+1} + I_2^{n+1} + I_3^{n+1} + I_4^{n+1} . 
\end{equation}
We proceed by bounding $I_1^i$, $I_2^i$, $I_3^i$, and $I_4^i$,
starting with $I_1^i$.

Restricting the error equation in \cref{lem:errorequation} for
$\boldsymbol{v}_h$ with general index $i$, we find the error equation
\begin{equation*}
  a_h (\boldsymbol{e}_u^{h, i} , \boldsymbol{v}_h) + b_h (\boldsymbol{v}_h, \boldsymbol{e}_{p_T}^{h,i}) = 0
  \quad \forall \boldsymbol{v}_h\in \boldsymbol{V}_h.
\end{equation*}
By \cref{eq:b_h-inf-sup}, the above equality,
\cref{eq:a_h-continuity}, the equivalence between $\tnorm{\cdot}_v$
and $\tnorm{\cdot}_{v'}$ \cite[eq. (5.5)]{Wells:2011}, and
\cref{eq:a_h-coercivity},
\begin{equation}
  \begin{split}
    C\tnorm{\boldsymbol{e}_{p_T}^{h,i}}_q
    &\le \sup_{\boldsymbol{0}\neq\boldsymbol{v}_h \in \boldsymbol{V}_h}
    \frac{ b_h(\boldsymbol{v}_h, \boldsymbol{e}_{p_T}^{h,i}) }{\tnorm{\boldsymbol{v}_h}_{v}}
    = \sup_{\boldsymbol{0}\neq\boldsymbol{v}_h \in \boldsymbol{V}_h}
    \frac{-a_h (\boldsymbol{e}_u^{h, i} , \boldsymbol{v}_h)}{\tnorm{\boldsymbol{v}_h}_{v}}
    \\
    &\le C\mu\tnorm{\boldsymbol{e}_u^{h,i}}_{v}
    \le C\mu^{1/2}\del[1]{a_h(\boldsymbol{e}_u^{h,i}, \boldsymbol{e}_u^{h,i})}^{1/2},
  \end{split}
\end{equation}
implying that
\begin{equation}
  \label{eq:ept-estimate}
  \mu^{-1/2} \tnorm{\boldsymbol{e}_{p_T}^{h,i}}_q \le C A_i.
\end{equation}
We may now bound $I_1^i$ using \cref{eq:h_b-continuity} and
\cref{eq:ept-estimate}:
\begin{equation}
  \label{eq:boundingI1i}
  I_1^i
  \le \Delta t|b_h(d_t\boldsymbol{e}_u^{I,i},\boldsymbol{e}_{p_T}^{h,i})|
  \le C \Delta t \tnorm{d_t\boldsymbol{e}_u^{I,i}}_v \tnorm{\boldsymbol{e}_{p_T}^{h,i}}_q
  \le C \Delta t \mu^{1/2} \tnorm{d_t\boldsymbol{e}_u^{I,i}}_v A_i.
\end{equation}
A bound for $I_2^i$ follows from the Cauchy--Schwarz and Young's
inequalities:
\begin{equation*}
  I_2^i
  \le \kappa^{-1}\Delta t\|e_z^{I,i}\|_{\Omega}\|e_z^{h,i}\|_{\Omega}
  \leq \frac{\kappa^{-1}}{2}\Delta t\|e_z^{I,i}\|^2_{\Omega}+ B_{i}^2.
\end{equation*}
Using the Cauchy--Schwarz and triangle inequalities, we bound $I_3^i$
as follows:
\begin{align*}
  I_3^i
  &\leq c_0\Delta t  \|\partial_tp^{i}-d_t\Pi_Qp^{i}\|_{\Omega}\|e_p^{h,i}\|_{\Omega}
  \\
  &\leq c_0 \Delta t (\|\partial_tp^{i}-d_tp^{i}\|_{\Omega}+\|d_t e_p^{I,i}\|_{\Omega})\|e_p^{h,i}\|_{\Omega}
  \\
  &\leq (2 c_0)^{1/2} \Delta t (\|\partial_tp^{i}-d_tp^{i}\|_{\Omega}+\|d_t e_p^{I,i}\|_{\Omega}) A_{i} .
\end{align*}
To estimate $I_4^i$, we first derive an auxiliary result. By the
assumption that $C_* \mu \le \lambda$ (see \cref{sec:biot}),
\begin{align*}
  \lambda^{-1} \alpha^2 \| e_p^{h, i} \|_{\Omega}^2
  &\le 2 \lambda^{-1} \del[2]{\|\alpha e_p^{h,i}-e_{p_T}^{h,i}\|_{\Omega}^2 }  + 2 \lambda^{-1} \| e_{p_T}^{h,i} \|_{\Omega}^2
  \\
  &\le 2 \lambda^{-1} \del[2]{\|\alpha e_p^{h,i}-e_{p_T}^{h,i}\|_{\Omega}^2 }  + 2 (C_* \mu)^{-1} \| e_{p_T}^{h,i} \|_{\Omega}^2.
\end{align*}
Combining this estimate with \cref{eq:ept-estimate} we obtain:
\begin{equation}
  \label{eq:lambdaalphaephibound}
  \lambda^{-1/2} \alpha \| e_p^{h,i} \|_{\Omega} \le C A_{i}.
\end{equation}
The Cauchy--Schwarz and triangle inequalities, together with
\cref{eq:lambdaalphaephibound} now imply
\begin{align*}
  I_4^i
  & \leq \lambda^{-1}\alpha\Delta t \|\alpha (\partial_tp^{i} - d_t \Pi_Q p^{i})
    -(\partial_tp_T^{i}-d_t\Pi_Q p_T^{i})\|_{\Omega}\|e_p^{h,i}\|_{\Omega}
  \\
  & \leq C \lambda^{-1/2} \Delta t\Big(\alpha\|\partial_tp^{i} - d_t p^{i}\|_{\Omega} + \alpha\|d_t e_p^{I,i}\|_{\Omega} +\|\partial_tp_T^{i}-d_t p_{T}^{i}\|_{\Omega} + \|d_t e_{p_T}^{I,i}\|_{\Omega}\Big) A_{i} .
\end{align*}
If we define 
\begin{equation*}
  \begin{split}
    E_{i} =
    & C \mu^{1/2}\Delta t \tnorm{d_t \boldsymbol{e}_u^{I,i}}_{v} + (2c_0)^{1/2} \Delta t ( \| \partial_tp^{i} - d_t p^{i} \|_{\Omega} + \| d_t e_p^{I,i} \|_{\Omega})
    \\
    & + C\lambda^{-1/2} \Delta t\Big(\alpha\|\partial_tp^{i} - d_t p^{i}\|_{\Omega} + \alpha\|d_t e_p^{I,i}\|_{\Omega} + \|\partial_tp_T^{i}-d_t p_{T}^{i}\|_{\Omega} + \|d_t e_{p_T}^{I,i}\|_{\Omega}\Big),
    \\
    D_{i}
    =& \frac{\kappa^{-1}}{2} \Delta t \| e_z^{I,i} \|_{\Omega}^2 ,    
  \end{split}
\end{equation*}
we find, using \cref{eq:An-Bn-ineq0}, that
\begin{equation}
  (A_{i+1}^2-A_i^2)+B_{i+1}^2 
  \le E_{i+1}A_{i+1}+D_{i+1}.
\end{equation}
Summing now for $i=0$ to $i=n-1$, and using $A_0 = 0$, we obtain,
after shifting indices and using $D_0 \ge 0$ and $B_0=0$,
\begin{equation}
  A_n^2 + \sum_{j=1}^{n} B_{j}^2 \le \sum_{j=1}^{n} E_{j}A_{j} + \sum_{j=1}^{n} D_{j}.
\end{equation}
Then, by \cref{lemma:not-gronwall} we obtain
\begin{equation}
  \label{eq:An-estimate}
  A_n + \Big(\sum_{i=1}^n B_i^2\Big)^{1/2}
  \leq  C \left( \sum_{i=1}^n E_i +  \Big(\sum_{i=1}^n D_i\Big)^{1/2}\right).  
\end{equation}
To prove \cref{eq:error-estimate1-a}, we therefore need to estimate
$\sum_{i=1}^n E_i$ and $(\sum_{i=1}^n D_i )^{1/2}$.

By \cref{eq:a_h-coercivity,eq:ell-interpolation-property}, we note that
\begin{equation*}
  \begin{split}
    \mu^{1/2}\Delta t\sum_{i=1}^n \tnorm{d_t \boldsymbol{e}_u^{I,i}}_v
    &\le C\Delta t\sum_{i=1}^n (a_h(d_t \boldsymbol{e}_u^{I,i}, d_t \boldsymbol{e}_u^{I,i}))^{1/2}
    \\
    &\le C\Delta t\mu^{1/2}\sum_{i=1}^n h^l\norm[0]{d_t u^i}_{l+1,\Omega}
    \le C \mu^{1/2} h^l\norm{\partial_tu}_{L^1(I;H^{l+1}(\Omega))},
  \end{split}
\end{equation*}
where we used $\Delta t d_t u^i = \int_{t^{i-1}}^{t^i} \partial_t u(s) ds$ for the last inequality. Using this estimate,
together with \cref{ineq:taylor3,ineq:taylor4}, we find:
{\fontsize{9.0}{10.5}\selectfont
  \begin{equation}
    \label{eq:Ei-sum-estimate}
    \begin{split}
      \sum_{i=1}^n E_{i} 
      =& C\sum_{i=1}^n\bigg[ \mu^{1/2} \Delta t \tnorm{d_t \boldsymbol{e}_u^{I,i}}_v + (2c_0)^{1/2}  \Delta t ( \| \partial_tp^{i} - d_t p^{i} \|_{\Omega} + \| d_t e_p^{I,i} \|_{\Omega})  
      \\
      & + \lambda^{-1/2} \Delta t\Big(\alpha\|\partial_tp^{i} - d_t p^{i}\|_{\Omega} + \alpha\|d_t e_p^{I,i}\|_{\Omega} + \|\partial_tp_T^{i}-d_t p_{T}^{i}\|_{\Omega} + \|d_t e_{p_T}^{I,i}\|_{\Omega}\Big)\bigg] , 
      \\       
      \le& C \mu^{1/2} h^l\norm{\partial_tu}_{L^1(I;H^{l+1}(\Omega))}
      \\
      &+ C \max\cbr[0]{c_0^{1/2}, \lambda^{-1/2}}\del[1]{\Delta t\norm[0]{\partial_{tt}p,\partial_{tt}p_T}_{L^1(I;L^2(\Omega))} + h^l\norm[0]{\partial_tp,\partial_tp_T}_{L^1(I;H^l(\Omega))}}.
     \end{split}
\end{equation}
}
Next, by \cref{eq:interpolation-property},
\begin{equation}
  \label{eq:D-sum-estimate}
  \begin{split}
    \sum_{i=1}^n D_i
    = \sum_{i=1}^n \frac{\kappa^{-1}}{2} \Delta t \| e_z^{I,i} \|_{\Omega}^2
    &\le C \kappa^{-1} (\Delta t) h^{2l} \sum_{i=1}^n \| z (t^i) \|_{H^l(\Omega)}^2
    \\
    &\le C \kappa^{-1} T h^{2l} \| z \|_{C^0(I; H^l(\Omega))}^2,
  \end{split}
\end{equation}
where in the last inequality $n\Delta t \le T$ is used. Combining \cref{eq:An-estimate} with
\cref{eq:Ei-sum-estimate,eq:D-sum-estimate} and the coercivity of
$a_h$ \cref{eq:a_h-coercivity}, we find:
\begin{equation}
  \label{eq:approxerror}
  \begin{split}
    c_0^{1/2}&\norm[0]{e_p^{h,n}}_{\Omega} + \lambda^{-1/2}\norm[0]{\alpha e_p^{h,n} - e_{p_T}^{h,n}}_{\Omega}
    + \mu^{1/2}\tnorm{\boldsymbol{e}_u^{h,n}}_v + \kappa^{-1/2} \del[2]{\sum_{i=1}^n\Delta t\norm[0]{e_z^{h,i}}_{\Omega}^2}^{1/2}
    \\
    \le&
    C \mu^{1/2} h^l\norm{\partial_tu}_{L^1(I;H^{l+1}(\Omega))}
    \\
    &+ C \max\cbr[0]{c_0^{1/2}, \lambda^{-1/2}}\del[1]{\Delta t\norm[0]{\partial_{tt}p,\partial_{tt}p_T}_{L^1(I;L^2(\Omega))} + h^l\norm[0]{\partial_tp,\partial_tp_T}_{L^1(I;H^l(\Omega))}}
    \\
    &+ C \kappa^{-1/2} T^{1/2} h^{l} \| z \|_{C^0(I; H^l(\Omega))}
    \\
    \le& c_1 \Delta t + c_2 h^l    
  \end{split}
\end{equation}
where
\begin{align*}
  c_1
  =& C\max\cbr[0]{c_0^{1/2}, \lambda^{-1/2}}\norm[0]{\partial_{tt}p,\partial_{tt}p_T}_{L^1(I;L^2(\Omega))}
  \\  
  \le& C\max\cbr[0]{c_0^{1/2}, \lambda^{-1/2}}\norm[0]{p,p_T}_{W^{2,1}(I;L^2(\Omega))},
  \\
  c_2
  =& C \big(\mu^{1/2} \norm{\partial_tu}_{L^1(I;H^{l+1}(\Omega))}
     + \max\cbr[0]{c_0^{1/2}, \lambda^{-1/2}}\norm[0]{\partial_tp,\partial_tp_T}_{L^1(I;H^l(\Omega))}
  \\
   &+ \kappa^{-1/2} T^{1/2} \| z \|_{C^0(I; H^l(\Omega))} \big).
\end{align*}
Next, by \cref{eq:errorsplit}, the triangle inequality, and the
approximation error \cref{eq:approxerror},
\begin{equation}
  \label{eq:errorfirsttry}
  \begin{split}
    &c_0^{1/2}\norm[0]{p^n-p_h^n}_{\Omega} + \lambda^{-1/2}\norm[0]{\alpha(p^n-p_h^n) - (p_T^n - p_{Th}^n)}_{\Omega}
    + \mu^{1/2}\tnorm{\boldsymbol{u}^n - \boldsymbol{u}_h^n}_v
    \\
    &+ \kappa^{-1/2} \del[2]{\sum_{i=1}^n\Delta t\norm[0]{z^i-z_h^i}_{\Omega}^2}^{1/2}
    \\
    \le
    &
    c_0^{1/2}\norm[0]{e_p^{I,n}}_{\Omega} + \lambda^{-1/2}\norm[0]{\alpha e_p^{I,n} - e_{p_T}^{I,n}}_{\Omega}
    + \mu^{1/2}\tnorm{\boldsymbol{e}_u^{I,n}}_v + \del[2]{2 \sum_{i=1}^n D_i}^{1/2}
    \\    
    &+ c_1\Delta t + c_2 h^l.
  \end{split}
\end{equation}
Note that
\begin{align*}
  c_0^{1/2}\norm[0]{e_p^{I,n}}_{\Omega}
  &\le C c_0^{1/2} h^l\norm[0]{p^n}_{l,\Omega},
  \\
  \lambda^{-1/2}\norm[0]{\alpha e_p^{I,n} - e_{p_T}^{I,n}}_{\Omega}
  &\le C \lambda^{-1/2} h^l (\norm[0]{p^n}_{l,\Omega} + \norm[0]{p_T^n}_{l,\Omega}),
  \\
  \mu^{1/2}\tnorm{\boldsymbol{e}_u^{I,n}}_v
  &\le C\mu^{1/2} h^l\norm[0]{u^n}_{l+1,\Omega}.
\end{align*}
When combined with \cref{eq:D-sum-estimate,eq:errorfirsttry},
\begin{equation}
  \label{eq:errorfirsttry2}
  \begin{split}
    &c_0^{1/2}\norm[0]{p^n-p_h^n}_{\Omega} + \lambda^{-1/2}\norm[0]{\alpha(p^n-p_h^n) - (p_T^n - p_{Th}^n)}_{\Omega}
    + \mu^{1/2}\tnorm{\boldsymbol{u}^n - \boldsymbol{u}_h^n}_v
    \\
    &+ \kappa^{-1/2} \del[2]{\sum_{i=1}^n\Delta t\norm[0]{z^i-z_h^i}_{\Omega}^2}^{1/2}
    \le
    c_1\Delta t + (c_2 + c_3) h^l,
  \end{split}
\end{equation}
where
{\fontsize{8.5}{9.0}\selectfont
\begin{equation}
  \begin{split}
    c_3
    =& 
    C \del[2]{c_0^{1/2}\norm[0]{p^n}_{l,\Omega} + \lambda^{-1/2} (\norm[0]{p^n}_{l,\Omega} + \norm[0]{p_T^n}_{l,\Omega})
      + \mu^{1/2} \norm[0]{u^n}_{l+1,\Omega} + \kappa^{-1/2}T^{1/2} \norm[0]{z}_{C^0(I;H^l(\Omega))}}  
    \\
    =&
    C \del[2]{(c_0^{1/2}+\lambda^{-1/2})\norm[0]{p^n}_{l,\Omega} + \lambda^{-1/2}\norm[0]{p_T^n}_{l,\Omega}
      + \mu^{1/2} \norm[0]{u^n}_{l+1,\Omega} + \kappa^{-1/2}T^{1/2} \norm[0]{z}_{C^0(I;H^l(\Omega))}} .    
  \end{split}
\end{equation}
}
Let us have a closer look at $c_2+c_3$. Using the
Sobolev embedding
$W^{s,q}(I; H^l(\Omega)) \hookrightarrow C^0(I; H^l(\Omega))$ for
$(s,q) = (1,1)$ and $(s,q) = (2,1)$,
\begin{equation*}
  \begin{split}
    c_2+c_3
    \le&
    C \big(\mu^{1/2} \cbr[1]{\norm{\partial_tu}_{L^1(I;H^{l+1}(\Omega))} + \norm[0]{u^n}_{l+1,\Omega}}
    + \kappa^{-1/2} T^{1/2} \| z \|_{C^0(I; H^l(\Omega))}
    \\
    &+ \max\cbr[0]{c_0^{1/2}, \lambda^{-1/2}}\cbr[1]{\norm[0]{\partial_tp,\partial_tp_T}_{L^1(I;H^l(\Omega))}
      + \norm[0]{p^n}_{l,\Omega} + \norm[0]{p_T^n}_{l,\Omega}} \big)
    \\
    \le&
    C \big(\mu^{1/2} \cbr[1]{\norm{\partial_tu}_{L^1(I;H^{l+1}(\Omega))} + \norm[0]{u}_{C^0(I; H^{l+1}(\Omega))}}
    + \kappa^{-1/2} T^{1/2} \| z \|_{C^0(I; H^l(\Omega))}
    \\
    &+ \max\cbr[0]{c_0^{1/2}, \lambda^{-1/2}}\cbr[1]{\norm[0]{\partial_tp,\partial_tp_T}_{L^1(I;H^l(\Omega))}
      + \norm[0]{p,p_T}_{C^0(I; H^l(\Omega))}} \big)
    \\
    \le&
    C \big(\mu^{1/2} \norm[0]{u}_{W^{1,1}(I; H^{l+1}(\Omega))}
    + \kappa^{-1/2} T^{1/2} \| z \|_{C^0(I; H^l(\Omega))}
    \\
    &+ \max\cbr[0]{c_0^{1/2}, \lambda^{-1/2}}\norm[0]{p,p_T}_{W^{1,1}(I; H^l(\Omega))}\big),
  \end{split}
\end{equation*}
proving \cref{eq:error-estimate1-a}.

Finally, \cref{eq:error-estimate1-b} follows from \cref{eq:errorsplit}
and the triangle inequality, \cref{eq:ept-estimate},
\cref{eq:approxerror}, and usual interpolation estimates for the
$L^2$-projection:
\begin{equation}
  \begin{split}
    \mu^{-1/2}\tnorm{\boldsymbol{p}_T^n - \boldsymbol{p}_{Th}^n}_{q}
    &\le
    \mu^{-1/2}\tnorm{\boldsymbol{e}_{p_T}^{I,n}}_q + C A_n
    \\
    &\le
    c_1\Delta t + \del[1]{C\mu^{-1/2}\norm[0]{p_T^n}_{l,\Omega} + c_2}h^l
    \\
    &\le
    c_1\Delta t + \del[1]{C\mu^{-1/2}\norm[0]{p_T}_{C^0(I,H^l(\Omega))} + c_2}h^l
  \end{split}
\end{equation}
Let us have a closer look at the constant in front of the second term:
\begin{equation*}
  \begin{split}
    C\mu^{-1/2} & \norm[0]{p_T}_{C^0(I,H^l(\Omega))} + c_2
    \\
    \le&
    C \big(\mu^{1/2} \norm{\partial_tu}_{L^1(I;H^{l+1}(\Omega))}
    + \max\cbr[0]{c_0^{1/2}, \lambda^{-1/2}}\norm[0]{\partial_tp,\partial_tp_T}_{L^1(I;H^l(\Omega))}
    \\
    &+ \kappa^{-1/2} T^{1/2} \| z \|_{C^0(I; H^l(\Omega))}
    + C\mu^{-1/2}\norm[0]{p_T}_{C^0(I,H^l(\Omega))}\big)
    \\
    \le&
    C \big(\mu^{1/2} \norm{\partial_tu}_{L^1(I;H^{l+1}(\Omega))}
    + \max\cbr[0]{c_0^{1/2}, \mu^{-1/2}}\norm[0]{p,p_T}_{W^{1,1}(I;H^l(\Omega))}
    \\
    &+ \kappa^{-1/2} T^{1/2} \| z \|_{C^0(I; H^l(\Omega))}\big),
  \end{split}
\end{equation*}
where in the last step we used that $C^*\lambda^{-1/2} \le
\mu^{1/2}$. This proves \cref{eq:error-estimate1-b}. 
\end{proof}

We end this section by noting that the estimates in \cref{thm:error-estimate} for the displacement, Darcy velocity, and total pressure are unconditionally robust in the incompressible limit $c_0 \to 0$ and $\lambda \to \infty$.

\section{Numerical examples}
\label{sec:numerical_examples}
We now validate our theoretical analysis. As stated previously in
\cref{rem:edghdg}, the analysis in this paper holds both for HDG and
EDG-HDG. As such, both methods are implemented using the
Netgen/NGSolve finite element library
\cite{Schoberl:1997,Schoberl:2014}. Numerical results are compared to
analytical solutions and some benchmark problems.

\subsection{Convergence rates for a static problem}
\label{ss:convergence_steadyproblem}
We consider a test case proposed in~\cite[Example
1]{Oyarzua:2016}. Consider the static Biot problem \cref{eq:biot} on
$\Omega$ where \cref{eq:biot_c} is replaced by
\begin{equation}
  \label{eq:static_biot}
  c_0  p + \lambda^{-1}\alpha (\alpha p - p_T) + \nabla\cdot z = g \quad \text{in } \Omega.
\end{equation}
We consider a domain $\Omega$ with four curved boundaries parametrized
as
\begin{equation*}
  \begin{split}
  \Gamma_1 &= \cbr[1]{\omega\in[0,1]:\ x_1 = \omega + \gamma\cos(\pi\omega)\sin(\pi\omega),\
    x_2 = -\gamma\cos(\pi\omega)\sin(\pi\omega)},
  \\
  \Gamma_2 &= \cbr[1]{\omega\in[0,1]:\ x_1 = 1 + \gamma\cos(\pi\omega)\sin(\pi\omega),\
    x_2 = \omega - \gamma\cos(\pi\omega)\sin(\pi\omega)},
  \\
  \Gamma_3 &= \cbr[1]{\omega\in[1,0]:\ x_1 = \omega + \gamma\cos(\pi\omega)\sin(\pi\omega),\
    x_2 = 1-\gamma\cos(\pi\omega)\sin(\pi\omega)},
  \\
  \Gamma_4 &= \cbr[1]{\omega\in[1,0]:\ x_1 = \gamma\cos(\pi\omega)\sin(\pi\omega),\
    x_2 = \omega-\gamma\cos(\pi\omega)\sin(\pi\omega)},
  \end{split}
\end{equation*}
with $\gamma=-0.08$. We then define
$\Gamma_D = \Gamma_1 \cup \Gamma_3 \cup \Gamma_4$,
$\Gamma_P = \Gamma_1 \cup \Gamma_2$, $\Gamma_T = \Gamma_2$, and
$\Gamma_F = \Gamma_3 \cup \Gamma_4$. The solution to the Biot problem
is taken as
\begin{equation}
  \label{eq:tc1_solution}
  u = a
  \begin{bmatrix}
    \sin(\pi x_1)\cos(\pi x_2) + x_1^2/(2\lambda)
    \\
    -\cos(\pi x_1)\sin(\pi x_2) + x_2^2/(2\lambda)
  \end{bmatrix},
  \quad
  p = b\sin(\pi x_1)\sin(\pi x_2).
\end{equation}
This solution~\cref{eq:tc1_solution} is used to set the body force
$f$, the source/sink term $g$, and \emph{inhomogeneous} boundary
conditions. As parameters we set $a=10^{-4}$, $b=\pi$,
$\kappa=10^{-7}$, $\alpha=0.1$ and $c_0=10^{-5}$. We consider both
mild incompressibility ($\nu=0.4$) and quasi-incompressibility
($\nu=0.49999$) and consider two values for $E$, namely $E=10^4$ and
$E=1$. Furthermore, we consider the rates of convergence for the
lowest order ($k=1$) and a higher order ($k=3$) approximation (with
$k$ the polynomial approximation in \cref{eq:approximationspaces}).

\Cref{thm:error-estimate} does not present an estimate in the
$L^2$-norm for the displacement and only a suboptimal estimate for the
Darcy velocity. Nevertheless, since our main objective here is to show
the robustness of the discretization in the incompressible limit, we
present in \cref{tab:rates_tc1_opt0,tab:rates_tc1_opt1}, for the HDG
and EDG-HDG schemes, respectively, the errors and rates of convergence
of all unknowns in the $L^2$-norm. Let us first observe that the
velocity and displacement converge at rate $k+1$ and that the
pressures converge at rate $k$. These are optimal rates of
convergence. We furthermore observe that the errors for all unknowns
are independent of the value of Poisson's ratio $\nu$ and for the
modulus of elasticity $E$. It is particularly interesting to note that
the choice $E=10^4$ and $\nu=0.49999$ (corresponding to
$\lambda \approx 1.7\cdot 10^8$) does not affect the quality of the
approximation. This confirms the robustness of the error estimates in \Cref{thm:error-estimate} in the incompressible limit.

\begin{table}[tbp]
  \caption{Rates of convergence for HDG the test case described
    in~\cref{ss:convergence_steadyproblem} for $E=1$ or $E=10^4$,
    $\nu=0.4$ or $\nu=0.49999$, and for $k=1$ or $k=3$. Here $r$ is the
    rate of convergence.
  }
  {\small
    \begin{center}
      \begin{tabular}{ccccccccc}
        \hline
        Cells & $\norm{u_h - u}_{\Omega}$ & $r$ & $\norm{p_{Th} - p_T}_{\Omega}$ & $r$ & $\norm{z_h - z}_{\Omega}$
        & $r$ & $\norm{p_h - p}_{\Omega}$ & $r$ \\
        \hline
        \multicolumn{9}{l}{$k=1$, $E=10^4$, $\nu=0.4$} \\
    384 & 4.2e-07 & 2.0 & 4.7e-02 & 1.1 & 5.6e-09 & 2.1 & 1.3e-01 & 1.0 \\
   1536 & 1.1e-07 & 2.0 & 2.2e-02 & 1.1 & 1.3e-09 & 2.1 & 6.3e-02 & 1.0 \\
   6144 & 2.6e-08 & 2.0 & 1.0e-02 & 1.1 & 3.3e-10 & 2.0 & 3.2e-02 & 1.0 \\
  24576 & 6.6e-09 & 2.0 & 5.0e-03 & 1.0 & 8.1e-11 & 2.0 & 1.6e-02 & 1.0 \\
        \multicolumn{9}{l}{$k=1$, $E=10^4$, $\nu=0.49999$} \\
    384 & 4.3e-07 & 2.0 & 6.7e-02 & 1.2 & 3.9e-09 & 2.1 & 1.3e-01 & 1.0 \\
   1536 & 1.1e-07 & 2.0 & 2.9e-02 & 1.2 & 9.6e-10 & 2.0 & 6.3e-02 & 1.0 \\
   6144 & 2.7e-08 & 2.0 & 1.3e-02 & 1.1 & 2.4e-10 & 2.0 & 3.2e-02 & 1.0 \\
  24576 & 6.7e-09 & 2.0 & 6.3e-03 & 1.1 & 5.9e-11 & 2.0 & 1.6e-02 & 1.0 \\
        \multicolumn{9}{l}{$k=3$, $E=10^4$, $\nu=0.4$} \\
    384 & 3.2e-10 & 4.0 & 9.1e-05 & 3.0 & 5.0e-12 & 4.0 & 1.8e-04 & 3.0 \\
   1536 & 2.0e-11 & 4.0 & 1.1e-05 & 3.1 & 3.0e-13 & 4.0 & 2.3e-05 & 3.0 \\
   6144 & 1.2e-12 & 4.0 & 1.3e-06 & 3.0 & 1.8e-14 & 4.0 & 2.8e-06 & 3.0 \\
  24576 & 7.6e-14 & 4.0 & 1.7e-07 & 3.0 & 1.1e-15 & 4.0 & 3.5e-07 & 3.0 \\
        \multicolumn{9}{l}{$k=3$, $E=10^4$, $\nu=0.49999$} \\
    384 & 3.6e-10 & 4.0 & 1.7e-04 & 3.1 & 2.6e-12 & 4.0 & 1.8e-04 & 3.0 \\
   1536 & 2.2e-11 & 4.0 & 2.0e-05 & 3.1 & 1.6e-13 & 4.0 & 2.3e-05 & 3.0 \\
   6144 & 1.4e-12 & 4.0 & 2.4e-06 & 3.1 & 1.0e-14 & 4.0 & 2.8e-06 & 3.0 \\
  24576 & 8.5e-14 & 4.0 & 2.9e-07 & 3.0 & 6.4e-16 & 4.0 & 3.5e-07 & 3.0 \\
        \hline
        \multicolumn{9}{l}{$k=1$, $E=1$, $\nu=0.4$} \\
    384 & 6.2e-07 & 3.6 & 1.3e-02 & 1.0 & 3.7e-09 & 2.1 & 1.3e-01 & 1.0 \\
   1536 & 1.1e-07 & 2.5 & 6.3e-03 & 1.0 & 9.6e-10 & 2.0 & 6.3e-02 & 1.0 \\
   6144 & 2.6e-08 & 2.1 & 3.2e-03 & 1.0 & 2.6e-10 & 1.9 & 3.2e-02 & 1.0 \\
  24576 & 6.6e-09 & 2.0 & 1.6e-03 & 1.0 & 7.2e-11 & 1.9 & 1.6e-02 & 1.0 \\
        \multicolumn{9}{l}{$k=1$, $E=1$, $\nu=0.49999$} \\
    384 & 6.4e-07 & 3.7 & 1.3e-02 & 1.0 & 3.9e-09 & 2.1 & 1.3e-01 & 1.0 \\
   1536 & 1.1e-07 & 2.5 & 6.3e-03 & 1.0 & 9.5e-10 & 2.0 & 6.3e-02 & 1.0 \\
   6144 & 2.7e-08 & 2.0 & 3.2e-03 & 1.0 & 2.4e-10 & 2.0 & 3.2e-02 & 1.0 \\
  24576 & 6.7e-09 & 2.0 & 1.6e-03 & 1.0 & 5.9e-11 & 2.0 & 1.6e-02 & 1.0 \\
        \multicolumn{9}{l}{$k=3$, $E=1$, $\nu=0.4$} \\
    384 & 3.4e-10 & 4.1 & 1.8e-05 & 3.0 & 2.2e-12 & 3.7 & 1.8e-04 & 3.0 \\
   1536 & 2.0e-11 & 4.1 & 2.3e-06 & 3.0 & 2.0e-13 & 3.5 & 2.3e-05 & 3.0 \\
   6144 & 1.2e-12 & 4.0 & 2.8e-07 & 3.0 & 1.6e-14 & 3.7 & 2.8e-06 & 3.0 \\
  24576 & 7.6e-14 & 4.0 & 3.5e-08 & 3.0 & 1.1e-15 & 3.9 & 3.5e-07 & 3.0 \\
        \multicolumn{9}{l}{$k=3$, $E=1$, $\nu=0.49999$} \\
    384 & 3.6e-10 & 4.0 & 1.8e-05 & 3.0 & 2.6e-12 & 4.0 & 1.8e-04 & 3.0 \\
   1536 & 2.2e-11 & 4.0 & 2.3e-06 & 3.0 & 1.6e-13 & 4.0 & 2.3e-05 & 3.0 \\
   6144 & 1.4e-12 & 4.0 & 2.8e-07 & 3.0 & 1.0e-14 & 4.0 & 2.8e-06 & 3.0 \\
  24576 & 8.5e-14 & 4.0 & 3.5e-08 & 3.0 & 6.4e-16 & 4.0 & 3.5e-07 & 3.0 \\
        \hline
      \end{tabular}
      \label{tab:rates_tc1_opt0}
    \end{center}
  }
\end{table}

\begin{table}[tbp]
  \caption{Rates of convergence for EDG-HDG the test case described
    in~\cref{ss:convergence_steadyproblem} for $E=1$ or $E=10^4$,
    $\nu=0.4$ or $\nu=0.49999$, and for $k=1$ or $k=3$. Here $r$ is the
    rate of convergence.}
  {\small
    \begin{center}
      \begin{tabular}{ccccccccc}
        \hline
        Cells & $\norm{u_h - u}_{\Omega}$ & $r$ & $\norm{p_{Th} - p_T}_{\Omega}$ & $r$ & $\norm{z_h - z}_{\Omega}$
        & $r$ & $\norm{p_h - p}_{\Omega}$ & $r$ \\
        \hline
        \multicolumn{9}{l}{$k=1$, $E=10^4$, $\nu=0.4$} \\
    384 & 5.4e-07 & 2.1 & 7.1e-02 & 1.3 & 7.7e-09 & 2.2 & 1.3e-01 & 1.0 \\ 
   1536 & 1.3e-07 & 2.0 & 2.9e-02 & 1.3 & 1.7e-09 & 2.2 & 6.3e-02 & 1.0 \\ 
   6144 & 3.2e-08 & 2.0 & 1.3e-02 & 1.2 & 3.7e-10 & 2.2 & 3.2e-02 & 1.0 \\ 
  24576 & 8.1e-09 & 2.0 & 5.7e-03 & 1.1 & 8.7e-11 & 2.1 & 1.6e-02 & 1.0 \\         
        \multicolumn{9}{l}{$k=1$, $E=10^4$, $\nu=0.49999$} \\
    384 & 5.4e-07 & 2.1 & 1.6e-01 & 1.4 & 3.9e-09 & 2.1 & 1.3e-01 & 1.0 \\ 
   1536 & 1.3e-07 & 2.0 & 6.0e-02 & 1.4 & 9.6e-10 & 2.0 & 6.3e-02 & 1.0 \\ 
   6144 & 3.3e-08 & 2.0 & 2.3e-02 & 1.4 & 2.4e-10 & 2.0 & 3.2e-02 & 1.0 \\ 
  24576 & 8.1e-09 & 2.0 & 9.5e-03 & 1.3 & 5.9e-11 & 2.0 & 1.6e-02 & 1.0 \\         
        \multicolumn{9}{l}{$k=3$, $E=10^4$, $\nu=0.4$} \\
    384 & 3.4e-10 & 4.0 & 9.8e-05 & 3.1 & 5.3e-12 & 4.0 & 1.8e-04 & 3.0 \\ 
   1536 & 2.1e-11 & 4.0 & 1.2e-05 & 3.1 & 3.2e-13 & 4.1 & 2.3e-05 & 3.0 \\ 
   6144 & 1.3e-12 & 4.0 & 1.4e-06 & 3.0 & 1.9e-14 & 4.0 & 2.8e-06 & 3.0 \\ 
  24576 & 8.1e-14 & 4.0 & 1.7e-07 & 3.0 & 1.2e-15 & 4.0 & 3.5e-07 & 3.0 \\        
        \multicolumn{9}{l}{$k=3$, $E=10^4$, $\nu=0.49999$} \\
    384 & 4.0e-10 & 4.0 & 2.0e-04 & 3.2 & 2.6e-12 & 4.0 & 1.8e-04 & 3.0 \\ 
   1536 & 2.4e-11 & 4.0 & 2.3e-05 & 3.1 & 1.6e-13 & 4.0 & 2.3e-05 & 3.0 \\ 
   6144 & 1.5e-12 & 4.0 & 2.6e-06 & 3.1 & 1.0e-14 & 4.0 & 2.8e-06 & 3.0 \\ 
  24576 & 9.2e-14 & 4.0 & 3.1e-07 & 3.1 & 6.4e-16 & 4.0 & 3.5e-07 & 3.0 \\         
        \hline
        \multicolumn{9}{l}{$k=1$, $E=1$, $\nu=0.4$} \\
    384 & 7.0e-07 & 3.3 & 1.3e-02 & 1.0 & 3.9e-09 & 2.0 & 1.3e-01 & 1.0 \\ 
   1536 & 1.3e-07 & 2.4 & 6.3e-03 & 1.0 & 1.1e-09 & 1.9 & 6.3e-02 & 1.0 \\ 
   6144 & 3.2e-08 & 2.1 & 3.2e-03 & 1.0 & 3.0e-10 & 1.8 & 3.2e-02 & 1.0 \\ 
  24576 & 8.0e-09 & 2.0 & 1.6e-03 & 1.0 & 7.8e-11 & 1.9 & 1.6e-02 & 1.0 \\        
        \multicolumn{9}{l}{$k=1$, $E=1$, $\nu=0.49999$} \\
    384 & 7.2e-07 & 3.4 & 1.3e-02 & 1.0 & 3.9e-09 & 2.1 & 1.3e-01 & 1.0 \\ 
   1536 & 1.4e-07 & 2.4 & 6.3e-03 & 1.0 & 9.5e-10 & 2.0 & 6.3e-02 & 1.0 \\ 
   6144 & 3.3e-08 & 2.1 & 3.2e-03 & 1.0 & 2.4e-10 & 2.0 & 3.2e-02 & 1.0 \\ 
  24576 & 8.1e-09 & 2.0 & 1.6e-03 & 1.0 & 5.9e-11 & 2.0 & 1.6e-02 & 1.0 \\        
        \multicolumn{9}{l}{$k=3$, $E=1$, $\nu=0.4$} \\
    384 & 3.7e-10 & 4.1 & 1.8e-05 & 3.0 & 2.3e-12 & 3.6 & 1.8e-04 & 3.0 \\ 
   1536 & 2.2e-11 & 4.1 & 2.3e-06 & 3.0 & 2.1e-13 & 3.5 & 2.3e-05 & 3.0 \\ 
   6144 & 1.3e-12 & 4.0 & 2.8e-07 & 3.0 & 1.6e-14 & 3.7 & 2.8e-06 & 3.0 \\ 
  24576 & 8.1e-14 & 4.0 & 3.5e-08 & 3.0 & 1.1e-15 & 3.9 & 3.5e-07 & 3.0 \\ 
        \multicolumn{9}{l}{$k=3$, $E=1$, $\nu=0.49999$} \\
    384 & 4.0e-10 & 4.0 & 1.8e-05 & 3.0 & 2.6e-12 & 4.0 & 1.8e-04 & 3.0 \\ 
   1536 & 2.4e-11 & 4.0 & 2.3e-06 & 3.0 & 1.6e-13 & 4.0 & 2.3e-05 & 3.0 \\ 
   6144 & 1.5e-12 & 4.0 & 2.8e-07 & 3.0 & 1.0e-14 & 4.0 & 2.8e-06 & 3.0 \\ 
  24576 & 9.2e-14 & 4.0 & 3.5e-08 & 3.0 & 6.4e-16 & 4.0 & 3.5e-07 & 3.0 \\         
        \hline
      \end{tabular}
      \label{tab:rates_tc1_opt1}
    \end{center}
  }
\end{table}

\subsection{Convergence rates for the quasi-static problem}
\label{ss:convergence_unsteadyproblem}
We now consider a manufactured solution for the quasi-static
problem~\cref{eq:biot} on the unit square. We divide the boundary of
our domain into
\begin{align*}
  \Gamma_1 &= \cbr[0]{(x_1,x_2) \in \partial\Omega\,:\, x_2=0},
  &
  \Gamma_2 &= \cbr[0]{(x_1,x_2) \in \partial\Omega\,:\, x_1=1},
  \\
  \Gamma_3 &= \cbr[0]{(x_1,x_2) \in \partial\Omega\,:\, x_2=1},
  &
  \Gamma_4 &= \cbr[0]{(x_1,x_2) \in \partial\Omega\,:\, x_1=0},  
\end{align*}
and set $\Gamma_D = \Gamma_1 \cup \Gamma_3 \cup \Gamma_4$,
$\Gamma_P = \Gamma_1 \cup \Gamma_2$, $\Gamma_T = \Gamma_2$, and
$\Gamma_F = \Gamma_3 \cup \Gamma_4$. As exact solution we take
\begin{equation}
  \label{eq:tc3_solution}
  u =
  \begin{bmatrix}
    \sin(\pi t)\sin(\pi x_1)\sin(\pi x_2)\\
    \sin(\pi t)\sin(\pi x_1)\cos(\pi x_2)
  \end{bmatrix},
  \quad
  p = \sin(\pi (x_1 - x_2 - t)),
\end{equation}
and set body force terms, source/sink terms, initial and boundary
conditions accordingly. As parameters we set $E=10^4$,
$\kappa=10^{-2}$, $\alpha=0.1$, $c_0=0.1$, $\nu=0.2$. We consider the
solution over the time interval $I=(0,0.1]$ and show the rates of
convergence at $t=0.1$ in \cref{tab:rates_tc3} for HDG and EDG-HDG
using $k=1$ and $k=2$. We are interested in the spatial rates of convergence and so we implement a second order backward differentiation formulae (BDF2) time stepping
scheme and take a time step of
$\Delta t = 10^{-3}$ so that spatial errors dominate over temporal errors. In \cref{tab:rates_tc3}, we observe optimal
rates of convergence for all unknowns, both for the HDG
and EDG-HDG schemes. Furthermore, note that although the error in the
total pressure is relatively large, this has no effect on the errors
in the displacement, velocity and pore pressure of the fluid, which
are all magnitudes smaller.
\begin{table}[tbp]
  \caption{Rates of convergence for HDG and EDG-HDG for the test case
    described in~\cref{ss:convergence_unsteadyproblem} for $k=1$ and
    $k=2$. Here dofs are the total number of degrees-of-freedom and
    $r$ is the rate of convergence.}
  {\small
    \begin{center}
      \begin{tabular}{ccccccccc}
        \hline
        Dofs & $\norm{u_h - u}_{\Omega}$ & $r$ & $\norm{p_{Th} - p_T}_{\Omega}$ & $r$ & $\norm{z_h-z}_{\Omega}$ & $r$ & $\norm{p_h-p}_{\Omega}$ & $r$ \\
        \hline
        \\
        \multicolumn{9}{c}{HDG} \\
        \\
        \hline
        \multicolumn{9}{l}{$k=1$} \\
   896 & 1.7e-02 & 2.5 & 7.2e+02 & 1.1 & 4.1e-03 & 0.9 & 1.9e-01 & 0.7 \\ 
  3456 & 4.1e-03 & 2.1 & 3.6e+02 & 1.0 & 1.1e-03 & 1.9 & 9.3e-02 & 1.0 \\ 
 13568 & 1.0e-03 & 2.0 & 1.8e+02 & 1.0 & 3.0e-04 & 1.9 & 4.6e-02 & 1.0 \\ 
 53760 & 2.5e-04 & 2.0 & 9.0e+01 & 1.0 & 7.6e-05 & 2.0 & 2.3e-02 & 1.0 \\ 
        \multicolumn{9}{l}{$k=2$} \\
  1632 & 1.7e-03 & 1.6 & 1.1e+02 & 1.0 & 3.5e-04 & 3.3 & 2.8e-02 & 2.3 \\ 
  6336 & 2.1e-04 & 3.0 & 2.7e+01 & 2.0 & 4.3e-05 & 3.0 & 7.0e-03 & 2.0 \\ 
 24960 & 2.7e-05 & 3.0 & 6.8e+00 & 2.0 & 5.3e-06 & 3.0 & 1.8e-03 & 2.0 \\ 
 99072 & 3.4e-06 & 3.0 & 1.7e+00 & 2.0 & 6.8e-07 & 3.0 & 4.4e-04 & 2.0 \\        
        \hline
        \\
        \multicolumn{9}{c}{EDG-HDG} \\
        \\
        \hline
        \multicolumn{9}{l}{$k=1$} \\
   722 & 2.2e-02 & 2.3 & 7.4e+02 & 1.1 & 4.9e-03 & 0.4 & 1.8e-01 & 0.6 \\ 
  2786 & 5.5e-03 & 2.0 & 3.7e+02 & 1.0 & 1.2e-03 & 2.0 & 9.2e-02 & 1.0 \\ 
 10946 & 1.3e-03 & 2.0 & 1.8e+02 & 1.0 & 3.0e-04 & 2.0 & 4.6e-02 & 1.0 \\ 
 43394 & 3.3e-04 & 2.0 & 9.2e+01 & 1.0 & 7.5e-05 & 2.0 & 2.3e-02 & 1.0 \\
        \multicolumn{9}{l}{$k=2$} \\
  1458 & 1.8e-03 & 1.6 & 1.1e+02 & 1.1 & 4.4e-04 & 3.0 & 2.8e-02 & 2.3 \\ 
  5666 & 2.4e-04 & 2.9 & 2.8e+01 & 2.0 & 5.5e-05 & 3.0 & 7.0e-03 & 2.0 \\ 
 22338 & 3.0e-05 & 3.0 & 6.9e+00 & 2.0 & 6.7e-06 & 3.0 & 1.8e-03 & 2.0 \\ 
 88706 & 3.7e-06 & 3.0 & 1.7e+00 & 2.0 & 8.5e-07 & 3.0 & 4.4e-04 & 2.0 \\ 
        \hline
      \end{tabular}
      \label{tab:rates_tc3}
    \end{center}
  }
\end{table}

\subsection{The footing problem}
\label{ss:footing_problem}
The two-dimensional footing problem has been proposed in the
literature to study the locking-free properties of numerical methods
for the Biot equations \cite{Gaspar:2008,Oyarzua:2016}. We follow here
the setup of \cite{Ambartsumyan:2020} and consider the domain
$\Omega = (-50, 50)\times (0,75)$ and model parameters
$\kappa=10^{-4}$, $c_0=10^{-3}$, $\alpha=0.1$, $E=3\cdot 10^4$, and
$\nu = 0.4995$ (so that $\lambda \approx 10^7$). We define the
boundaries
$\Gamma_1=\cbr{(x_1,x_2)\in\partial\Omega,\ |x_1|\le 50/3, x_2=75}$,
$\Gamma_2=\cbr{(x_1,x_2)\in\partial\Omega,\ |x_1| > 50/3, x_2=75}$,
and $\Gamma_3 = \partial\Omega\backslash(\Gamma_1\cup\Gamma_2)$ and
impose the following boundary conditions:
\begin{equation*}
  \sigma n = (0, -\sigma_0)^T\ \text{on}\ \Gamma_1,\quad
  \sigma n = 0\ \text{on}\ \Gamma_2,\quad
  u =0\ \text{on}\ \Gamma_3, \quad
  p = 0\ \text{on}\ \partial\Omega,
\end{equation*}
where $\sigma_0 = 10^4$. As initial conditions we impose $u(x,0) = 0$
and $p(x,0) = 0$. We solve this problem with HDG until $T=50$ using
BDF2 time stepping. We choose a time step of size $\Delta t = 1$, take
$k=2$ in our finite element spaces, and compute the solution on an
unstructured mesh consisting of 169984 simplices.

\begin{figure}
  \centering
  \subfloat[Displacement $u$. \label{fig:footing_u}]{\includegraphics[width=0.47\textwidth]{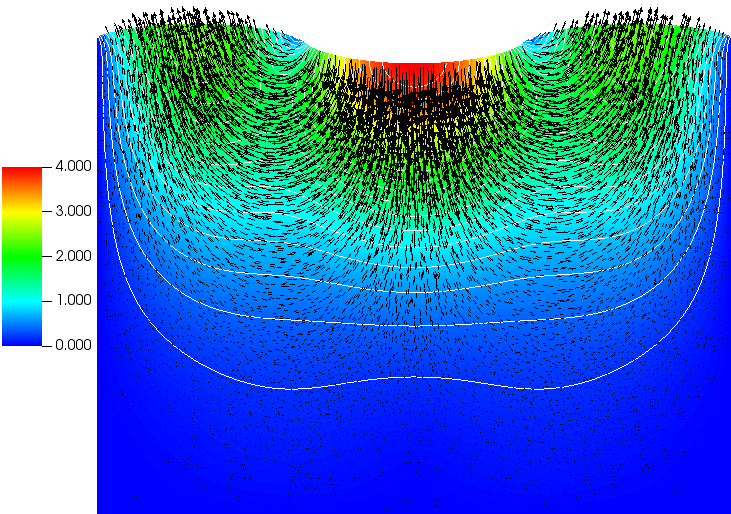}}
  \
  \subfloat[Darcy velocity $z$. \label{fig:footing_z}]{\includegraphics[width=0.49\textwidth]{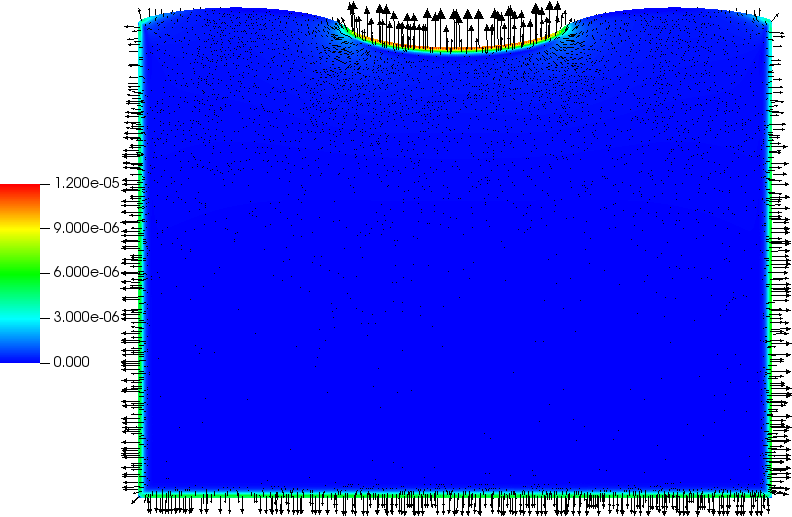}}
  \\
  \subfloat[Fluid pressure $p$. \label{fig:footing_p}]{\includegraphics[width=0.48\textwidth]{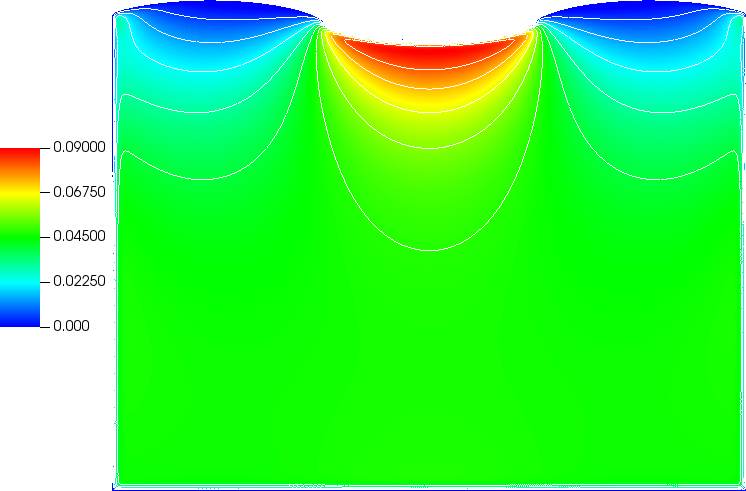}}
  \
  \subfloat[Total pressure $p_T$. \label{fig:footing_pT}]{\includegraphics[width=0.48\textwidth]{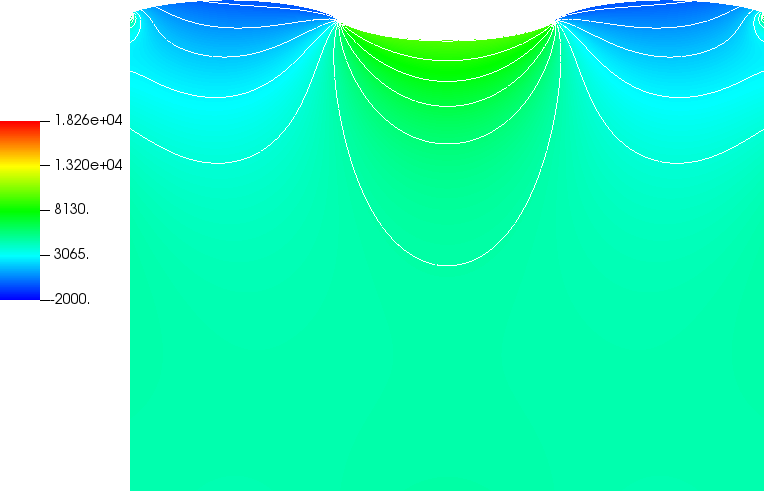}}
  \caption{The solution to the footing problem of \cref{ss:footing_problem} in the deformed domain at $t=50$. }
  \label{fig:plots_footing}
\end{figure}

We show the solution to this problem at time $t=50$ in
\cref{fig:plots_footing}. In this incompressible limit we observe that
the discretization results in pressure and displacement solutions are free of, respectively, spurious oscillations and locking
effects.

\subsection{The cantilever bracket problem}
\label{ss:cantilever_problem}
The cantilever bracket problem was used in
\cite{Ambartsumyan:2020,Liu:2004,Phillips:2009} to study locking
phenomena at low permeability and when the storage coefficient is
zero. Consider the domain $\Omega = (0,1)^2$ and define
\begin{align*}
  \Gamma_1 &= \cbr{(x_1,x_2) \in \partial\Omega,\ x_2=0},
  &
    \Gamma_2 &= \cbr{(x_1,x_2) \in \partial\Omega,\ x_1=1},
  \\
  \Gamma_3 &= \cbr{(x_1,x_2) \in \partial\Omega,\ x_2=1},
  &
  \Gamma_4 &= \cbr{(x_1,x_2) \in \partial\Omega,\ x_1=0}.
\end{align*}
We impose the boundary conditions
\begin{equation*}
  z\cdot n = 0 \text{ on } \partial\Omega,\quad
  \sigma n = (0, -1)^T \text{ on } \Gamma_3,\quad
  \sigma n = 0 \text{ on } \Gamma_1 \cup \Gamma_2,\quad
  u = 0 \text{ on } \Gamma_4.
\end{equation*}
At $t=0$ we set $u=0$ and $p=0$. The model parameters are chosen as
$E=10^5$, $\nu=0.4$, $\alpha = 0.93$, $c_0=0$, and $\kappa = 10^{-7}$
\cite{Phillips:2009}. As shown in \cite{Phillips:2009}, with these
parameters continuous Galerkin numerical methods show spurious
oscillations in the pressure on a very short time
interval. Therefore, we consider here the time interval
$I=[0,0.005]$. In our discretization we combine the EDG-HDG
discretization with BDF2 time stepping, choose a time step of
$\Delta t = 0.001$, set $k=2$ in our finite element spaces, and
compute the solution on a mesh consisting of 128 simplices.

\begin{figure}
  \centering
  \subfloat[Pressure $p$ at $t=0.001$. \label{fig:cantilever-p}]{\includegraphics[width=0.47\textwidth]{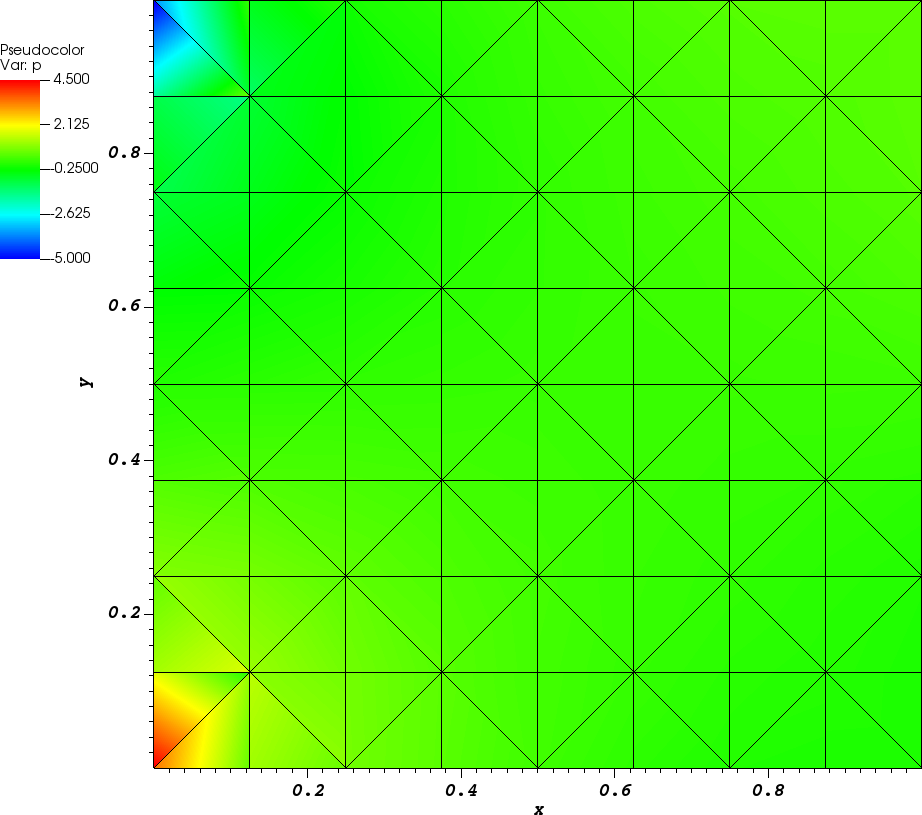}}
  \
  \subfloat[Pressure $p$ along different $x$-lines. \label{fig:cantilever-x}]{\includegraphics[width=0.49\textwidth]{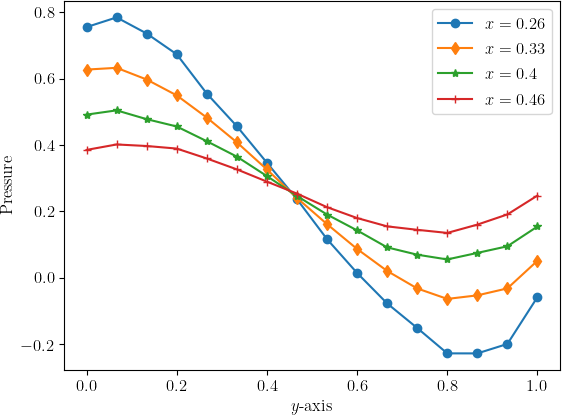}}
  \caption{The solution to the cantilever bracket problem
    \cref{ss:cantilever_problem} using EDG-HDG. Left: the pressure
    solution at $t=0.001$. Right: the pressure solution along
    different $x$-lines at time $t=0.005$. }
  \label{fig:cantilever}
\end{figure}

We plot the solution in \cref{fig:cantilever}. In
\cref{fig:cantilever-p} we observe that the pressure field at
$t=0.001$ is free from spurious oscillations, similar to the
discontinuous Galerkin solutions obtained in \cite{Phillips:2009}. We
further show in \cref{fig:cantilever-x} that the pressure solution
along the lines $x=0.26$, $x=0.33$, $x=0.4$, and $x=0.46$ at $t=0.005$
is free of oscillations, agreeing with other stable finite element
methods for this problem
\cite{Ambartsumyan:2020,Liu:2004,Phillips:2009}.

\section{Conclusions}
\label{sec:conclusions}

An HDG and an EDG-HDG method have been presented and analyzed for the total pressure formulation of the quasi-static poroelasticity model. Both discretization methods are shown to be well-posed and space-time a priori error estimates show robustness of the proposed methods when $\lambda \to \infty$ and $c_0 \to 0$; both methods are free of volumetric locking. Numerical examples confirm our theory and further show optimal spatial rates of convergence in the $L^2$-norm.

\section*{Statements and Declarations}

\noindent {\bf Funding} For AC and JL this material is based upon work supported by the
National Science Foundation under grant numbers DMS-2110782 and
DMS-2110781. SR gratefully acknowledges support from the Natural
Sciences and Engineering Research Council of Canada through the
Discovery Grant program (RGPIN-05606-2015).

\bibliographystyle{amsplain}
\bibliography{references}
\section{Appendix} 
\subsection{The inf-sup condition for $b_h$}
\label{ap:inf-sup_bhwh0qh}

By definition of $b_h$ \cref{eq:b_h},
\begin{equation}
  b_h((w_h,0), \boldsymbol{q}_h) :=
  \underbrace{- (q_h, \nabla\cdot w_h )_{\Omega}}_{=:b_h^1(w_h,q_h)}
  + \underbrace{\langle \bar{q}_h, w_h \cdot n\rangle_{\partial \mathcal{T}}}_{=:b_h^2(w_h,\bar{q}_h)}
  \quad \forall w_h\in V_h, \forall \boldsymbol{q}_h \in \boldsymbol{Q}_h^0.
\end{equation}
Let $q_h\in Q_h$. It is known (see, for example, \cite[Section
4.1.4]{Ern:book} or \cite[Remark 3.3]{JJLee:2017}) that there exists
$w \in \sbr[0]{H^1_{0,\Gamma_F}(\Omega)}^{d} := \cbr[0]{v \in
  \sbr[0]{H^1(\Omega)}^d \, :\, v|_{\Gamma_F}=0}$ such that
\begin{equation}
  \label{eq:inf-sup-bw}
  -(\nabla \cdot w, q_h)_{\Omega} =  \|q_h\|_{\Omega}^2, \quad \|w\|_{1,\Omega}\leq C\|q_h\|_{\Omega},
\end{equation}
for some positive constant $C$ that only depends on $\Omega$. Let
$\Pi_V: \sbr[0]{H^1(\Omega)}^d \rightarrow V_h$ be the BDM
interpolation operator \cite[Section III.3]{Brezzi:book} and observe
that by the single-valuedness of $\bar{q}_h$, continuity of
$\Pi_V w \cdot n$ across interior faces, and since $\bar{q}_h=0$ on
$\Gamma_P$ and $w=0$ on $\Gamma_F$,
\begin{equation*}
  b_h^2(\Pi_V w,\bar{q}_h) = \langle \bar{q}_h, \Pi_V w\cdot n\rangle_{\Gamma_F} = 0,  
\end{equation*}
i.e.,
$\Pi_V w\in \text{Ker}\, b_2 := \cbr[0]{ w_h\in V_h \, : \,
  b_h^2(w_h,\bar{q}_h) = 0 \quad \forall \bar{q}_h \in
  \bar{Q}_h}$. Recall also that
$(q_h,\nabla \cdot \Pi_V w)_{\Omega} = (q_h,\nabla \cdot w)_{\Omega}$
and $\|\Pi_V w\|_{\Omega} \leq C\|w\|_{1,\Omega}$. Then, by
\cref{eq:inf-sup-bw},
\begin{align*}
    \sup_{0\neq w_h \in  \text{Ker}\, b_2} \frac{b_h^1(w_h, q_h)}{\|w_h\|_{\Omega}}
    \geq \frac{-(q_h,\nabla \cdot \Pi_V w)_{\Omega}}{\|\Pi_V w\|_{\Omega}}
    \geq \frac{\|q_h\|^2_{\Omega}}{C\|q_h\|_{\Omega}}= C\|q_h\|_{\Omega}.
\end{align*}
Next, let $w_h:=L\bar{q}_h \in P_k(K)^{d}$ where $L$ is the local BDM
interpolation operator \cite{Brezzi:book} such that
\begin{equation}
    (L\bar{q}_h)\cdot n = \bar{q}_h, \quad \|L\bar{q}_h\|_K\leq Ch_K^{1/2}\|\bar{q}_h\|_{\partial K}, \quad K\in \mathcal{T}.
\end{equation}
Then 
\begin{align*}
  \sup_{0\neq w_h \in  V_h} \frac{b_h^2(w_h, \bar{q}_h)}{\|w_h\|_{\Omega}}
  & \geq \frac{\|\bar{q}_h\|^2_{\partial \mathcal{T}}}{\|w_h\|_{\Omega}}
    \geq \frac{\|\bar{q}_h\|^2_{\partial \mathcal{T}}}{C\sum_{K\in \mathcal{T}}h_K^{-1}\|w_h\|_{\Omega}}
  \\
  &\geq Ch_{\min}h_{\max}^{-1}\Big(\sum_{K\in \mathcal{T}}h_K\|\bar{q}_h\|_{\partial K}^2\Big)^{1/2}.
\end{align*}
Therefore, by \cite[Theorem 3.1]{Howell:2011},
\begin{align*}
  \sup_{0\neq w_h \in  V_h} \frac{b_h((w_h,0), \boldsymbol{q}_h)}{\|w_h\|_{\Omega}}
  &=\sup_{0\neq w_h \in  V_h} \frac{b_h^1(w_h, q_h)+b_h^2(w_h, \bar{q}_h)}{\|w_h\|_{\Omega}}\\
  &  \geq C\del[2]{ \|q_h\|_{\Omega}+h_{\min}h_{\max}^{-1}\del[1]{\sum_{K\in \mathcal{T}}h_K\|\bar{q}_h\|_{\partial K}^2}}
    \geq C\tnorm{\boldsymbol{q}_h}_q.
\end{align*}

\subsection{Error estimates following from Taylor's theorem}
\label{ap:err-taylor}

We prove here \cref{ineq:taylor3,ineq:taylor4}. Let
$D \subset \mathbb{R}^d$. Then for $\psi$ a regular enough function
defined on $[0,T] \times D$, using Taylor's theorem,
\begin{align*}
  \Delta t  \|\partial_t\psi^m-d_t\psi^m\|_{0,D} & =\| \Delta t \partial_t\psi^m - (\psi^{m} - \psi^{m-1})\|_{0,D} \\
  &= \| \int_{t_{m-1}}^{t_m} \partial_{tt}\psi(t)\underbrace{(t-t_{m-1})}_{\le \Delta t} \dif t \|_{0,D} \\
  &\le \Delta t \int_{t_{m-1}}^{t_m} \| \partial_{tt}\psi(t) \|_{0,D} \dif t \\
  &= \Delta t \| \partial_{tt}\psi \|_{L^1(t_{m-1},t_m;L^2(D))},
\end{align*}
from which \cref{ineq:taylor3} follows. Next, to show
\cref{ineq:taylor4} we use the identity
\begin{equation*}
  \Delta t d_t e_{\psi}^{I,m} = \int_{t^{m-1}}^{t^{m}} (\partial_t\psi(s) - \Pi_Q \partial_t\psi(s)) \dif s,
\end{equation*}
and \cref{eq:interpolation-property}. Then, by the approximation
property of the $L^2$-projection,
\begin{equation*}
  \begin{split}
    \Delta t \sum_{m=1}^n \| d_t e_{\psi}^{I,m} \|_{0,D}
    &\le  \sum_{m=1}^n \int_{t^{m-1}}^{t^{m}} \| \partial_t\psi(s) - \Pi_Q \partial_t\psi(s)) \|_{0,D} \dif s
    \\
    &\le C h^l \| \partial_t\psi \|_{L^1(I; H^l(D))}.
  \end{split}
\end{equation*}

\providecommand{\bysame}{\leavevmode\hbox to3em{\hrulefill}\thinspace}
\providecommand{\MR}{\relax\ifhmode\unskip\space\fi MR }
\providecommand{\MRhref}[2]{%
  \href{http://www.ams.org/mathscinet-getitem?mr=#1}{#2}
}
\providecommand{\href}[2]{#2}

\end{document}